\documentclass[12pt,twoside,leqno]{article}
\usepackage{amsmath}
\usepackage{amssymb}
\usepackage{amsxtra}
\usepackage{amscd}
\usepackage{amsthm}
\usepackage[mathscr]{eucal}

\theoremstyle{plain}

\newtheorem{sbthm}[subsubsection]{Theorem}
\newtheorem{sbprop}[subsubsection]{Proposition}
\newtheorem{sbcor}[subsubsection]{Corollary}
\newtheorem{sblem}[subsubsection]{Lemma}

\theoremstyle{definition}

\newtheorem{para}[subsection]{}

\newtheorem{sbrem}[subsubsection]{Remark}
\newtheorem{sbconj}[subsubsection]{Conjecture}
\newtheorem{sbpara}[subsubsection]{}

\newenvironment{pf}{\proof[\proofname]}{\endproof}

\begin{document}

\title{On log motives}

\author{Tetsushi Ito, Kazuya Kato, Chikara Nakayama, Sampei Usui}
\date{}
\maketitle
\newcommand\Cal{\mathcal}
\newcommand\define{\newcommand}

\define\gp{\mathrm{gp}}%
\define\fs{\mathrm{fs}}%
\define\an{\mathrm{an}}%
\define\mult{\mathrm{mult}}%
\define\Ker{\mathrm{Ker}\,}%
\define\Coker{\mathrm{Coker}\,}%
\define\Hom{\mathrm{Hom}\,}%
\define\Ext{\mathrm{Ext}\,}%
\define\rank{\mathrm{rank}\,}%
\define\gr{\mathrm{gr}}%
\define\cHom{\Cal{Hom}}
\define\cExt{\Cal Ext\,}%

\define\cA{\Cal A}
\define\cB{\Cal B}
\define\cC{\Cal C}
\define\cD{\Cal D}
\define\cO{\Cal O}
\define\cS{\Cal S}
\define\cM{\Cal M}
\define\cG{\Cal G}
\define\cH{\Cal H}
\define\cE{\Cal E}
\define\cF{\Cal F}
\define\cN{\Cal N}
\define\cT{\Cal T}
\define\fF{\frak F}
\define\Dc{\check{D}}
\define\Ec{\check{E}}
\define\cR{\Cal R}

\newcommand{\N}{{\mathbb{N}}}
\newcommand{\Q}{{\mathbb{Q}}}
\newcommand{\Z}{{\mathbb{Z}}}
\newcommand{\R}{{\mathbb{R}}}
\newcommand{\C}{{\mathbb{C}}}
\newcommand{\bN}{{\mathbb{N}}}
\newcommand{\bQ}{{\mathbb{Q}}}
\newcommand{\bF}{{\mathbb{F}}}
\newcommand{\bZ}{{\mathbb{Z}}}
\newcommand{\bP}{{\mathbb{P}}}
\newcommand{\bR}{{\mathbb{R}}}
\newcommand{\bC}{{\mathbb{C}}}
\newcommand{\bbQ}{{\bar \mathbb{Q}}}
\newcommand{\ol}[1]{\overline{#1}}
\newcommand{\too}{\longrightarrow}
\newcommand{\respect}{\rightsquigarrow}
\newcommand{\compatible}{\leftrightsquigarrow}
\newcommand{\upc}[1]{\overset {\lower 0.3ex \hbox{${\;}_{\circ}$}}{#1}}
\newcommand{\Gmlog}{\bG_{m, \log}}%
\newcommand{\Gm}{\bG_m}%
\newcommand{\ep}{\varepsilon}
\newcommand{\Spec}{\operatorname{Spec}}
\newcommand{\val}{{\mathrm{val}}} 
\newcommand{\n}{\operatorname{naive}}
\newcommand{\bs}{\operatorname{\backslash}}
\newcommand{\Gal}{\operatorname{{Gal}}}
\newcommand{\gal}{{\rm {Gal}}({\bar \Q}/{\Q})}
\newcommand{\galp}{{\rm {Gal}}({\bar \Q}_p/{\Q}_p)}
\newcommand{\gall}{{\rm{Gal}}({\bar \Q}_\ell/\Q_\ell)}
\newcommand{\wep}{W({\bar \Q}_p/\Q_p)}
\newcommand{\wel}{W({\bar \Q}_\ell/\Q_\ell)}
\newcommand{\Ad}{{\rm{Ad}}}
\newcommand{\BS}{{\rm {BS}}}
\newcommand{\even}{\operatorname{even}}
\newcommand{\End}{{\rm {End}}}
\newcommand{\odd}{\operatorname{odd}}
\newcommand{\GL}{\operatorname{GL}}
\newcommand{\np}{\text{non-$p$}}
\newcommand{\g}{{\gamma}}
\newcommand{\G}{{\Gamma}}
\newcommand{\Lam}{{\Lambda}}
\newcommand{\La}{{\Lambda}}
\newcommand{\lam}{{\lambda}}
\newcommand{\la}{{\lambda}}
\newcommand{\uL}{{{\hat {L}}^{\rm {ur}}}}
\newcommand{\uQp}{{{\hat \Q}_p}^{\text{ur}}}
\newcommand{\sel}{\operatorname{Sel}}
\newcommand{\dt}{{\rm{Det}}}
\newcommand{\Sig}{\Sigma}
\newcommand{\fil}{{\rm{fil}}}
\newcommand{\SL}{{\rm{SL}}}
\newcommand{\spl}{{\rm{spl}}}
\newcommand{\st}{{\rm{st}}}
\newcommand{\Isom}{{\rm {Isom}}}
\newcommand{\Mor}{{\rm {Mor}}}
\newcommand{\bg}{\bar{g}}
\newcommand{\id}{{\rm {id}}}
\newcommand{\cone}{{\rm {cone}}}
\newcommand{\al}{a}
\newcommand{\ChL}{{\cal{C}}(\La)}
\newcommand{\Image}{{\rm {Image}}}
\newcommand{\toric}{{\operatorname{toric}}}
\newcommand{\torus}{{\operatorname{torus}}}
\newcommand{\Aut}{{\rm {Aut}}}
\newcommand{\Qp}{{\mathbb{Q}}_p}
\newcommand{\barQp}{{\mathbb{Q}}_p}
\newcommand{\Qpur}{{\mathbb{Q}}_p^{\rm {ur}}}
\newcommand{\Zp}{{\mathbb{Z}}_p}
\newcommand{\Zl}{{\mathbb{Z}}_l}
\newcommand{\Ql}{{\mathbb{Q}}_l}
\newcommand{\Qlur}{{\mathbb{Q}}_l^{\rm {ur}}}
\newcommand{\F}{{\mathbb{F}}}
\newcommand{\eps}{{\epsilon}}
\newcommand{\epsLa}{{\epsilon}_{\La}}
\newcommand{\epsLaVxi}{{\epsilon}_{\La}(V, \xi)}
\newcommand{\epsOLaVxi}{{\epsilon}_{0,\La}(V, \xi)}
\newcommand{\Qplin}{{\mathbb{Q}}_p(\mu_{l^{\infty}})}
\newcommand{\otimesQplin}{\otimes_{\Qp}{\mathbb{Q}}_p(\mu_{l^{\infty}})}
\newcommand{\galFl}{{\rm{Gal}}({\bar {\Bbb F}}_\ell/{\Bbb F}_\ell)}
\newcommand{\gallur}{{\rm{Gal}}({\bar \Q}_\ell/\Q_\ell^{\rm {ur}})}
\newcommand{\galFF}{{\rm {Gal}}(F_{\infty}/F)}
\newcommand{\galFv}{{\rm {Gal}}(\bar{F}_v/F_v)}
\newcommand{\galF}{{\rm {Gal}}(\bar{F}/F)}
\newcommand{\epsVxi}{{\epsilon}(V, \xi)}
\newcommand{\epsOVxi}{{\epsilon}_0(V, \xi)}
\newcommand{\plim}{\lim_
{\scriptstyle 
\longleftarrow \atop \scriptstyle n}}
\newcommand{\sig}{{\sigma}}
\newcommand{\ga}{{\gamma}}
\newcommand{\del}{{\delta}}
\newcommand{\Vss}{V^{\rm {ss}}}
\newcommand{\Bst}{B_{\rm {st}}}
\newcommand{\Dpst}{D_{\rm {pst}}}
\newcommand{\Dcrys}{D_{\rm {crys}}}
\newcommand{\DdR}{D_{\rm {dR}}}
\newcommand{\Fin}{F_{\infty}}
\newcommand{\Kla}{K_{\lambda}}
\newcommand{\Ola}{O_{\lambda}}
\newcommand{\Mla}{M_{\lambda}}
\newcommand{\Det}{{\rm{Det}}}
\newcommand{\Sym}{{\rm{Sym}}}
\newcommand{\LaSa}{{\La_{S^*}}}
\newcommand{\cX}{{\cal {X}}}
\newcommand{\MHG}{{\frak {M}}_H(G)}
\newcommand{\tauMla}{\tau(M_{\lambda})}
\newcommand{\Fvur}{{F_v^{\rm {ur}}}}
\newcommand{\Lie}{{\rm {Lie}}}
\newcommand{\cL}{{\cal {L}}}
\newcommand{\cW}{{\cal {W}}}
\newcommand{\fq}{{\frak {q}}}
\newcommand{\cont}{{\rm {cont}}}
\newcommand{\SC}{{SC}}
\newcommand{\Om}{{\Omega}}
\newcommand{\dR}{{\rm {dR}}}
\newcommand{\crys}{{\rm {crys}}}
\newcommand{\hatSig}{{\hat{\Sigma}}}
\newcommand{\rdet}{{{\rm {det}}}}
\newcommand{\ord}{{{\rm {ord}}}}
\newcommand{\BdR}{{B_{\rm {dR}}}}
\newcommand{\BdRO}{{B^0_{\rm {dR}}}}
\newcommand{\Bcrys}{{B_{\rm {crys}}}}
\newcommand{\Qw}{{\mathbb{Q}}_w}
\newcommand{\barkappa}{{\bar{\kappa}}}
\newcommand{\cP}{{\Cal {P}}}
\newcommand{\cZ}{{\Cal {Z}}}
\newcommand{\oppLa}{{\Lambda^{\circ}}}
\newcommand{\syn}{{{\rm syn}}}
\newcommand{\Tr}{{{\rm Tr}}}
\newcommand{\ch}{{{\rm char}}}
\newcommand{\bG}{{{\bold G}}}
\newcommand{\nnum}{{{\rm num}}}
\newcommand{\coh}{{{\rm coh}}}
\newcommand{\CH}{{{\rm CH}}}

\begin{abstract}
\noindent 
We define the categories of log motives and log mixed motives.  The latter gives a new formulation for the category of mixed motives.  We prove that the former is a semisimple abelian category if and only if the numerical equivalence and homological equivalence coincide, and that it is also equivalent to that the latter is a Tannakian category.  We discuss various realizations, formulate Tate and Hodge conjectures, and verify them in curve case.
\end{abstract}

\renewcommand{\thefootnote}{\fnsymbol{footnote}}

\medskip 

\noindent {\bf Contents}

\smallskip
\S\ref{intro}. Introduction

\S\ref{preparation}. Preparations on log geometry

\S\ref{s:logmot}. Log motives

\S\ref{s:logmm}. Log mixed motives

\S\ref{realization}. Formulation with various realizations

\S\ref{example}. Examples

\medskip

\section{Introduction}
\label{intro}
\begin{para}
In this paper, we define

\medskip

(1) the category of log motives over an fs log scheme, and

(2) the category of log mixed motives over an fs log scheme.

\medskip

(1) is a generalization of the category of Grothendieck motives over a
field with respect to the homological equivalence. The category (2)
has $\oplus, \otimes$, dual, kernel and cokernel. We prove that the
following (i), (ii), and (iii) are equivalent.

\medskip

(i) The numerical equivalence and homological equivalence coincide in
the category (1).

(ii) The category (1) is a semisimple abelian category.

(iii) The category (2) is a Tannakian category.

\medskip

The equivalence of (i) and (ii) is the log version of the famous theorem of Jannsen (\cite{Jannsen}). 
\end{para}

\begin{para}\label{1.1}  
  We explain the organization briefly. 
  In this paper, except in \ref{logcoh}, an fs log scheme means an fs log scheme which has charts Zariski locally. 

 Let $S$ be an fs log scheme. We fix a prime number $\ell$ and assume that $\ell$ is invertible over $S$. 

After we give preparations in Section \ref{preparation}, we define in Section \ref{s:logmot} the category of log motives over $S$, which is the log version of the category of motives of Grothendieck. 
 In Section \ref{s:logmm}, we define the category of log mixed motives over $S$ basing on the theory in Section \ref{s:logmot}. 

Here we work modulo homological equivalence using $\ell$-adic log \'etale cohomology theory. 

  Vologodsky also defined log motives (cf.\ http://www.ihes.fr/\,$\tilde{}$\,abbes/Ogus/ogus70-programme.html).

  In the case where the log structure of $S$ is trivial, our construction gives a category of mixed motives over $S$ modulo homological equivalence. 
This does not use the theory of Voevodsky (\cite{Voe}), though we hope our theory is connected to it. 
In the case $S=\Spec(k)$ for a field $k$ of characteristic $0$ with trivial log structure, our definition of the category of mixed motives over $S$ is different from the definition of the category of mixed motives over $k$ given by Jannsen in \cite{Jannsen0}.
  The difference lies in that in the definition of morphisms, we use $K$-theory whereas he uses absolute Hodge cycles. 

  In Section \ref{realization}, we introduce realizations other than $\ell$-adic one. 
  In Section \ref{example}, we discuss examples. 
\end{para}

\begin{para}\label{1.2} We explain each section of this paper more. 

  We explain more about Section \ref{preparation}. 
  In Section \ref{preparation}, we give preparations on log geometry. We review results on log \'etale cohomology, log Betti cohomology, 
log de Rham cohomology, and log Hodge theory in \ref{logcoh}, and then review or prove  results on fans (\ref{fan}), on log modifications (\ref{logmodification}), and on the Grothendieck group of vector bundles on log schemes 
(\ref{Grothendieckgp}). 
\end{para}

\begin{para}\label{1.3}  We explain more about Section \ref{s:logmot}. 

Fix a prime number $\ell$ and let $S$ be an fs log scheme on which $\ell$ is invertible. 
We define the category of log motives over $S$ by imitating the definition of motive by Grothendieck modulo homological equivalence.

Recall that for a field $k$ whose characteristic is not $\ell$, the category of motives over $k$ modulo ($\ell$-adic) homological equivalence is defined as follows (cf.\ \cite{Scholl}). 
  For a projective smooth scheme $X$ over $k$ and for $r\in \Z$, consider a symbol  $h(X)(r)$. For projective smooth schemes $X$, $Y$ over $k$ and for $r,s\in\Z$, by a morphism $h(X)(r) \to h(Y)(s)$, we mean a homomorphism $\bigoplus_i \; H^i(X)_{\ell}(r) \to \bigoplus_i \; H^i(Y)_{\ell}(s)$ which comes from $\CH(X \times Y)_{\Q}$. Here $H^i(X)_{\ell}$ is the \'etale cohomology group $H^i_{\mathrm{\acute{e}t}}(X\otimes_k \bar k, \Q_{\ell})$ with $\bar k$ a fixed separable closure of $k$, $(r)$ denotes the $r$-th Tate twist, the same for $Y$ and $s$, and where $\CH=\bigoplus_i\;  \CH^i$ is the Chow group and $(\;)_{\Q}$ means $\otimes \Q$. A motive over $k$ is a pair $(h(X)(r), e)$, where $X$ is a projective smooth scheme over $k$, $r\in \Z$,  and $e$ is an idempotent of the ring of endomorphisms of $h(X)(r)$.

Imitating this, we define the category of log motives over $S$ is as follows. (See \ref{LM} for details.)
For a projective vertical log smooth fs log scheme $X$ over $S$ and for $r\in \Z$, consider a symbol $h(X)(r)$. For projective vertical log smooth fs log schemes $X$, $Y$ over $S$ and for $r,s\in \Z$, by a morphism $h(X)(r) \to h(Y)(s)$, we mean a homomorphism $h: \bigoplus_i \; H^i(X)_{\ell}(r) \to \bigoplus_i \; H^i(Y)_{\ell}(s)$ satisfying the condition (C) below. Here $H^i(X)_{\ell}$ is the smooth $\Q_{\ell}$-sheaf on the log \'etale site on $S$ defined to be the $i$-th relative log \'etale cohomology of $X$ over $S$,   $(r)$ denotes the $r$-th Tate twist, and the same for $Y$ and $s$.

\medskip

(C) For any geometric standard log point $p$ (\ref{slp}) over $S$, the pullback of $h$ to $p$ comes from an element of $\bigoplus_i\;  \gr^iK(Z)_{\Q}$ for  some log modification $Z$ of $X_p \times_p Y_p$, where $K(Z)$ denotes the Grothendieck group of the category of vector bundles on $Z$ and $\gr^i$ denotes the $i$-th graded quotient for the $\gamma$-filtration (\cite{SGA6}).

\medskip
A log motive over $S$ is a pair $(h(X)(r), e)$, where $X$  over $S$ and $r$ are as above and $e$ is an idempotent of the ring of endomorphisms of 
$h(X)(r)$ (\ref{deflmot}). 

The reason why we need log modifications is explained in \ref{iddia2}. 

In the case where $S=\Spec(k)$ for a field $k$ with the trivial log structure, by the fact $\gr^i K(Z)_\Q= \CH^i(Z)_\Q$ for any  smooth scheme $Z$ over $k$,  we have that our category of log motives over $S$ coincides with the category of motives over $k$ modulo homological equivalence due to Grothendieck.

\medskip

We will also define the category of log motives over $S$ modulo numerical equivalence by taking the quotient of the set of morphisms by numerical equivalence. We prove the following log version of the theorem of Jannsen.

\medskip

\noindent 
{\bf Theorem (= Theorem \ref{thmss})}. (1) The category of log motives over $S$ modulo numerical equivalence is a semisimple abelian category. 

\medskip

(2) The category of log motives over $S$ (defined in \ref{LM}) is a semisimple abelian category if and only if the numerical equivalence for morphisms of this category is trivial. 
\end{para}

\begin{para}
\label{1.4} 
We explain more about Section \ref{s:logmm}. Let $S$ and $\ell$ be as in \ref{1.3}. Roughly speaking, we follow the method of Deligne (\cite{D1}, \cite{D2}), who constructed mixed Hodge structures of geometric origin by using only projective smooth schemes over $\C$. 

Our definition of log mixed motives  is rather simple  and  is easily obtained 
by using the category of log (pure) motives in Section \ref{s:logmot}. 
  This may seem strange because usually it is impossible to take care of mixed objects by using only pure objects. 
  The reason why such a simple definition is fine is explained in \ref{justification}.

We will prove the following result.

\medskip

\noindent 
{\bf Theorem (cf.\ Theorem \ref{mthm}).}
Assume that the category of log motives over $S$ is semisimple, that is, the numerical equivalence coincides with the homological 
equivalence for this category (see (2) of the previous theorem). 
  Then the category of log mixed motives over $S$ is a  Tannakian category.
  In particular, it is an abelian category. 

\end{para}

\begin{para}\label{1.5}
In Sections \ref{preparation}--\ref{s:logmm}, our discussion only uses $\ell$-adic \'etale realization. 
  We consider in Section \ref{realization} more realizations, and formulate Tate conjecture and Hodge conjecture for log mixed motives.
  In the last section \ref{example}, we prove that these conjectures are true in certain cases (Proposition \ref{H2}, Proposition \ref{propE1}, Proposition \ref{propE2}). To prove the results Proposition \ref{propE1} and Proposition \ref{propE2} on morphisms between $H^1$ of log curves,  we use the theory of log abelian varieties in \cite{KKN2} and the theory of log Jacobian varieties \cite{Kajiwara}. 
\end{para}

\noindent
{\bf Acknowledgments.}
A part of this paper was written while the first author was staying
at Institut des Hautes \'Etudes Scientifiques.
He would like to thank IHES for its hospitality.
  We thank Takeshi Saito and Takeshi Kajiwara for helpful discussions. 
  The second author is partially  supported by NSF grants DMS 1303421 and DMS 1601861.   
  The third author is partially supported by JSPS, Kakenhi (C) No.\ 22540011, 
(B) No.\ 23340008, and (C) No.\ 16K05093.
  The fourth author is partially supported by JSPS, Kakenhi
(B) No.\ 23340008 and (C) No.\ 17K05200.

\section{Preparations on log geometry}
\label{preparation}

Basic references on log geometry are
\cite{KatoLog}, \cite{IllusieLog}.
Basic references on log \'etale cohomology are
 \cite{Nakayama1}, \cite{Nakayama2}, \cite{IllusieOverview}.
Basic references on algebraic cycles and $K$-groups are
\cite{SGA6}, \cite{Fulton}.

In this paper, except in \ref{logcoh},  we consider fs log schemes which have charts Zariski locally. 
A {\it monoid} means a commutative semigroup with a unit element which is usually denoted by $1$. 

  Let $X$ be an fs log scheme over an fs log scheme $S$.  
  We say that $X$ is {\it projective} if the underlying scheme of $X$ is projective over the underlying scheme of $S$. 
  We say that $X$ is {\it vertical} if for any point $x$ of $X$, whose image in $S$ is denoted by $s$, the 
face of $M_{X,\overline x}$ spanned by the image of $M_{S,\overline s}$ is the whole $M_{X,\overline x}$. 
  Cf.\ \cite{Nakayama1} Definition and Notation (7.3). 

  A morphism $f\colon X \to Y$ of integral log schemes is {\it exact} if for any $x \in X$, an element of $M^{\gp}_{Y,\overline{f(x)}}$ whose image in 
$M^{\gp}_{X,\overline{x}}$ belongs to $M_{X,\overline{x}}$ belongs to $M_{Y,\overline{f(x)}}$.
  Cf.\ \cite{KatoLog} {\it Definition} (4.6). 
 
\subsection{Log cohomology theories}
\label{logcoh}

  We review some theorems on log \'etale cohomology, log Betti cohomology, and log de Rham cohomology. 

  First we discuss the theorems on log \'etale cohomology. 
  There are two kinds of log \'etale cohomologies, that is, the kummer log \'etale (k\'et) one and the full log \'etale one. 
  In this paper, we mainly use the full log \'etale (l\'et) one.

  Let $f\colon X \to S$  be a morphism of fs log schemes. 
  Let $\ell$ be a prime number which is invertible on $S$. 
  Let $\Lambda = \bZ/\ell^n \bZ$ $(n\geq 1)$. 

\begin{sbprop}\label{prsm}
  Assume that $f\colon X \to S$ is proper and log smooth. 
  Then $R^q f_{\mathrm{l\acute{e}t}*}\Lambda$ (the higher direct image for the full log \'etale topology) is locally constant 
and constructible (see {\rm{\cite{Nakayama3}}}, $8.1$ for the definition) for all $q \in \bZ$. 
\end{sbprop}

\begin{pf}
  This is by \cite{Nakayama3}, Theorem 13.1 (1). 
\end{pf}

\begin{sbpara}
\label{smoothQ_l}
  As in the classical case, we define a {\it constructible 
$\bZ_{\ell}$-sheaf} as an inverse system $(F_n)_{n}$, where 
$F_n$ is a constructible sheaf of $\bZ/\ell^{n+1}\bZ$-modules such that 
$\bZ/\ell^n\bZ \otimes F_{n} \overset \cong \to F_{n-1}$. 
  A {\it smooth $\bZ_{\ell}$-sheaf} is a constructible $\bZ_{\ell}$-sheaf 
$(F_n)_n$ with each $F_n$ locally constant. 
  The smooth $\bZ_{\ell}$-sheaves make an abelian category. 
  We define the category of {\it constructible $\bQ_{\ell}$-sheaves} as the localization 
  of this abelian category 
by torsion sheaves.  
  By the above proposition, we have, under the assumption there, 
a smooth $\Q_{\ell}$-sheaf on $S_{\mathrm{l\acute{e}t}}$, which we denote by 
$R^q f_{\mathrm{l\acute{e}t}*}\bQ_{\ell}$.
\end{sbpara}

  The next is the Poincar\'e duality. 

\begin{sbprop}\label{Poin} 
  Let $d\geq0$. 
  Assume that $f\colon X \to S$ is proper, log smooth, vertical, and, full log \'etale locally on $S$, 
all fibers are of equi-$d$-dimensional.
  Then there is a natural isomorphism 
$$R^{2d-i}f_{\mathrm{l\acute{e}t}*}\Lambda(d)\overset{\cong} \to {\mathcal H}om\,(R^if_{\mathrm{l\acute{e}t}*}\Lambda, \Lambda)$$ 
for any $i$. 
\end{sbprop}

\begin{pf}
  This is by \cite{Nakayama3}, Theorem 14.2 (3). 
\end{pf}

\begin{sbcor}\label{Q_lPoin} 
  Under the same assumptions, suppose further that $S$ is noetherian. 
  Then, there is a natural isomorphism 
$$R^{2d-i}f_{\mathrm{l\acute{e}t}*}\bQ_{\ell}(d)\overset{\cong} \to {\mathcal H}om\,(R^if_{\mathrm{l\acute{e}t}*}\bQ_{\ell}, \bQ_{\ell})$$ 
for any $i$. 
\end{sbcor}

  The third is the K\"unneth formula. 

\begin{sbprop}\label{Kunn}
  Assume that $S$ is quasi-compact and that $f\colon X \to S$ is proper.
  Let $g:Y\to S$ be another proper morphism of fs log schemes.
  Let $h$ be the induced morphism $X \times_SY \to S$.
  Then there is a natural isomorphism
$$Rf_{\mathrm{l\acute et}*}\Lambda \otimes_{\Lambda}^{\bold L}Rg_{\mathrm{l\acute et}*}\Lambda \overset \cong \to
Rh_{\mathrm{l\acute et}*}\Lambda.$$
\end{sbprop}

\begin{pf}
  This is by \cite{Nakayama3}, Theorem 9.1.
\end{pf}

As a corollary, we have 

\begin{sbcor}\label{Q_lKunn}
  Assume that $S$ is quasi-compact and that $f\colon X \to S$ is proper and log smooth.
  Let $g:Y\to S$ be another proper and log smooth morphism of fs log schemes.
  Let $h$ be the induced morphism $X \times_SY \to S$. 
  Then, for each $n\geq0$, there is a natural isomorphism 
$$\bigoplus_{p+q=n}
R^pf_{\mathrm{l\acute et}*}\bQ_{\ell}\otimes
R^qg_{\mathrm{l\acute et}*}\bQ_{\ell}\overset \cong \to 
R^nh_{\mathrm{l\acute et}*}\bQ_{\ell}.$$
\end{sbcor}

\begin{pf}
  The natural homomorphism is seen to be bijective at stalks by 
the previous proposition.  
\end{pf}

  Next the theorems on log Betti cohomology are as follows. 
  Let $f\colon X \to S$  be a morphism of fs log analytic spaces.

\begin{sbprop}\label{prsmB}
  Assume that $f\colon X \to S$ is proper (i.e., the underlying map is universally closed and separated) and log smooth. 
  Then $R^q f^{\log}_{*}\bZ$ is locally constant sheaf of finitely generated abelian groups for all $q \in \bZ$. 
\end{sbprop}

\begin{pf}
  This is \cite{KajiwaraNakayama}, {\sc Corollary} 0.3.
\end{pf}

\begin{sbprop}\label{PoinB}
Let $d\geq0$. 
  Assume that $f\colon X \to S$ is proper, log smooth, vertical, and 
all fibers are of equi-$d$-dimensional.
  Then there is a natural isomorphism 
$$R^{2d-i}f^{\log}_{*}\bQ\overset{\cong} \to {\mathcal H}om\,(R^if^{\log}_{*}\bQ, \bQ)$$ 
for any $i$. 
\end{sbprop}

\begin{pf}
  The case where $f$ is exact is by \cite{NakayamaOgus}, Theorem 5.10 (3).
  The general case is reduced to this case by exactification as follows. 
  First, we assume that $S$ has a chart by an fs monoid and fix such a chart. 
  Then, by exactification (\cite{IllusieKatoNakayama} {\sc Proposition} (A.4.4)), there is a log blow-up (\cite{IllusieKatoNakayama} {\sc Definition} (6.1.1)) $p\colon S' \to S$ such that 
the base-changed morphism $f'\colon X':=X\times_SS' \to S'$ is exact. 
  By the exact case, we have the natural isomorphism 

\smallskip

\noindent $(*)$ $\qquad\qquad R^{2d-i}f^{\prime\log}_{*}\bQ\overset{\cong} \to {\mathcal H}om\,(R^if^{\prime\log}_{*}\bQ, \bQ)$

\smallskip

\noindent 
on $S'{}^{\log}$.  
  Below we will prove that sending this by $p^{\log}_*$ gives us an isomorphism 
$R^{2d-i}f^{\log}_{*}\bQ\overset{\cong} \to {\mathcal H}om\,(R^if^{\log}_{*}\bQ, \bQ)$ 
on $S^{\log}$.
  To see that the last isomorphism is independent of the choices of log blow-ups, we can argue as in \cite{Nakayama3} (14.10), where 
the $\ell$-adic analog of the same problem is treated. 
  Then, it implies that the isomorphism is independent also of the choices of charts, and glues into 
the desired isomorphism. 

  Now we calculate $p^{\log}_*$ of each side of $(*)$. 
  Since $R^{j}f^{\prime\log}_{*}\bQ$ is locally constant for any $j$ (Proposition \ref{prsmB}), 
by \cite{KajiwaraNakayama} {\sc Proposition} 5.3 (2), we have 
$p^{\log}_*Rf^{\prime\log}_{*}\bQ=Rp^{\log}_*Rf^{\prime\log}_{*}\bQ=Rf^{\log}_*p^{\log}_{*}\bQ=Rf^{\log}_*\bQ$, 
where we denote the base-changed morphism of $p$ by the same symbol and the last equality is by 
\cite{KajiwaraNakayama} {\sc Proposition 5.3} (1). 
  Hence, $p^{\log}_*R^{2d-i}f^{\prime\log}_{*}\bQ = R^{2d-i}f^{\log}_{*}\bQ$. 

  On the other hand, as for the right-hand-side of $(*)$, again by \cite{KajiwaraNakayama} {\sc Proposition} 5.3 (2), 
we have $R^if^{\prime \log}_*\bQ = p^{\log-1}p^{\log}_*R^if^{\prime\log}_*\bQ$, and it is isomorphic to 
$p^{\log-1}R^if^{\log}_*\bQ$ by the same argument for the left-hand-side. 
  Then, $p^{\log}_*{\mathcal H}om\,(R^if^{\prime\log}_{*}\bQ, \bQ)=
p^{\log}_*{\mathcal H}om\,(p^{\log-1}R^if^{\log}_{*}\bQ, \bQ)
={\mathcal H}om\,(R^if^{\log}_{*}\bQ, p^{\log}_*\bQ)
={\mathcal H}om\,(R^if^{\log}_{*}\bQ, \bQ)$, 
where the last equality is again by \cite{KajiwaraNakayama} {\sc Proposition} 5.3 (1).
  Thus we have an isomorphism 
$R^{2d-i}f^{\log}_{*}\bQ\overset{\cong} \to {\mathcal H}om\,(R^if^{\log}_{*}\bQ, \bQ)$. 
\end{pf}

\begin{sbprop}
\label{prsmbcB}
  Let $f\colon X \to S$ be a proper and log smooth morphism of fs log analytic spaces. 
  Let $g\colon S' \to S$ be any morphism of fs log analytic spaces. 
  Let $f'\colon X':=X \times_SS' \to S'$ and $g' \colon X' \to X$ be the base-changed morphisms.
  Let $L$ be a locally constant sheaf of abelian groups on $X^{\log}$. 
  Then the base change homomorphism $$g^{\log-1}Rf^{\log}_*L \to Rf^{\prime\log}_*g^{\log*}L$$ 
is an isomorphism. 
\end{sbprop}

\begin{pf}
  We may assume that $S$ has a chart. 
  By exactification (\cite{IllusieKatoNakayama} {\sc Proposition} (A.4.4)), we take a log blow-up $p\colon S_1 \to S$ such that 
the base-changed morphism $f_1\colon X_1:=X\times_SS_1 \to S_1$ is exact. 
  Then, by proper log smooth base change theorem in log Betti cohomology (\cite{KajiwaraNakayama} {\sc Theorem} 0.1), 
the cohomologies of $Rf_{1*}^{\log}p_X^{\log-1}L$ are locally constant, where $p_X$ is the base-changed morphism 
$X_1 \to X$. 
  Hence, by the invariance of cohomology under log blow-up (\cite{KajiwaraNakayama} {\sc Proposition} 5.3), 
to prove Proposition \ref{prsmbcB}, we can replace $f$ and $g$ by the base-changed ones with respect to $p$, and $L$ by 
its pullback $p_X^{\log-1}L$. 
  Thus we may assume that $f$ is exact. 
  Then the conclusion follows from the log proper base change theorem \cite{KajiwaraNakayama} 
{\sc Proposition} 5.1 (cf.\ \cite{KajiwaraNakayama} {\sc Remark} 5.1.1).
\end{pf}

\begin{sbprop}
\label{KunnB}
  Let the notation and assumption be as in the previous proposition.
  Assume that $g$ is proper. 
  Let $h\colon X' \to S$ be the induced morphism. 
  Then there is a natural isomorphism
$$Rf^{\log}_{\mathrm{an}*}\bQ \otimes_{\bQ}^{\bold L}Rg^{\log}_{\mathrm{an}*}\bQ \overset \cong \to
Rh^{\log}_{\mathrm{an}*}\bQ.$$
\end{sbprop}

\begin{pf}
  This is by Proposition \ref{prsmbcB} and the usual projection formula. 
\end{pf}

  Next is a comparison between log Betti cohomology and log \'etale cohomology. 

\begin{sbpara}
\label{slp} 
  A {\it standard log point} means the fs log scheme 
 $\Spec(k)$ for a field $k$ endowed with the log structure  associated to $\N\to k\;;\;1\mapsto 0$. If we like to present $k$, we call it  
a {\it standard log point associated to $k$}. 
The standard log point associated to an algebraically closed field is called a {\it geometric standard log point}.
\end{sbpara}

\begin{sbprop}
\label{comparison}
  Let $f\colon X \to S$ be a proper, log smooth and vertical morphism of fs log schemes 
with $S$ being of finite type over $\bC$. 
  Let 
$$X^{\log}_{\an} \overset \eta \to X_{\mathrm{k\acute{e}t}} \overset \kappa \gets 
X_{\mathrm{l\acute{e}t}}$$
be natural morphisms of sites (for $\eta$, see {\rm \cite{KatoNakayama}} Remark ($2.7$)). 
  Let $n \geq 1$ and $\Lambda = \bZ/\ell^n\bZ$.  
  Then we have 
$$\eta^*Rf_{\mathrm{k\acute{e}t}*}\Lambda = Rf^{\log}_{\an*}\Lambda, \quad
\kappa^*Rf_{\mathrm{k\acute{e}t}*}\Lambda = Rf_{\mathrm{l\acute{e}t}*}\Lambda.$$
\end{sbprop}

\begin{pf}
  The second one is shown in 13.4 of \cite{Nakayama3}.
  We prove the first one. 
  First, note that the cohomologies of the left-hand-side are locally constant and 
constructible by \cite{Nakayama3} Theorem 13.1 (2) and those of the right-hand-side are 
locally constant by Proposition \ref{prsmB}. 

  We reduce to the case where $f$ is exact. 
  We may assume that $S$ has a chart by an fs monoid and fix such a chart. 
  Then, by \cite{Nakayama3} Lemma 3.10, there is a log blow-up $p\colon S' \to S$ such that 
the base-changed morphism $f'\colon X':=X\times_SS' \to S'$ is exact. 
  By \cite{Nakayama3} Theorem 5.5 (1) and \cite{Nakayama3} Theorem 5.8 (1), 
we have $p^*_{\mathrm{k\acute{e}t}}Rf_{\mathrm{k\acute{e}t}*}\Lambda 
= p^*_{\mathrm{k\acute{e}t}}Rf_{\mathrm{k\acute{e}t}*}Rp_{\mathrm{k\acute{e}t}*}\Lambda 
= p^*_{\mathrm{k\acute{e}t}}Rp_{\mathrm{k\acute{e}t}*}Rf'_{\mathrm{k\acute{e}t}*}\Lambda 
= Rf'_{\mathrm{k\acute{e}t}*}\Lambda$, where we denote the base-changed morphism of $p$ by 
the same symbol. 

  Similarly, by \cite{KajiwaraNakayama} {\sc Proposition} 5.3, we have 
$p^{\log*}Rf^{\log}_{*}\Lambda 
= p^{\log*}Rf^{\log}_{*}Rp^{\log}_{*}\Lambda 
= p^{\log*}Rp^{\log}_{*}Rf^{\prime\log}_{*}\Lambda 
= Rf^{\prime\log}_{*}\Lambda$.
Thus we may and will assume that $f$ is exact. 

  Since the cohomologies of both sides are locally constant, we can work at stalks. 
  Let $s_0$ be a point of $S$. 
  By the following proposition \ref{desing}, there are a morphism $s \to S$ from the standard log point $s$ over $\bC$ 
whose image is $s_0$, and a log blow-up $X'$ of $X_s := X\times_Ss$ such that 
the composition $X' \to X_s \to s$ is strict semistable, i.e., a log deformation with smooth irreducible components. 
  It is enough to show that the homomorphism at a stalk over a point of $s_0^{\log}$ is bijective. 
  Then by the exact proper base change theorem (\cite{Nakayama1} Theorem (5.1) and Remark (5.1.1) 
for the log \'etale cohomology, \cite{KajiwaraNakayama} {\sc Proposition} 5.1, {\sc Remark} 5.1.1, and the 
usual proper base change theorem for topological spaces for the log Betti cohomology), we may assume that $S=s$, and further, 
by \cite{KajiwaraNakayama} {\sc Proposition} 5.3 (1) and \cite{Nakayama3} Theorem 5.5 (1), we may assume that 
$X=X'$, that is, in the original setting, we may assume that $S$ is the standard log point over $\bC$ and $X$ is 
strict semistable over $S$.

  Here we use the Steenbrink--Rapoport--Zink (SRZ, for short) spectral sequences as follows. 
  In the proof of \cite{FujisawaNakayama} Theorem 7.1, it is shown that 
there is a natural isomorphism 
between the $\ell$-adic SRZ spectral sequence and the Betti SRZ spectral sequence. 
  Since these converge to the stalk of $\ell$-adic log \'etale cohomologies and 
that of log Betti cohomologies, respectively, we have the desired isomorphism. 
\end{pf}

\begin{sbprop}
\label{desing}
  Let $s=(\Spec k, \bN)$ be a standard log point.
  Let $X \to s$ be a quasi-compact, vertical, and log smooth morphism of fs log schemes. 
  Then there are a positive integer $n$ and a log blow-up ({\rm \cite{Nakayama3}} $2.2$) $X' \to X \times_s{s_n}$, where 
$s_n:=(\Spec k, \frac 1{n} \bN)$, such that the composition $X' \to s_n$ is strict semistable. 
\end{sbprop}

  This is a variant of the semistable reduction theorem of D.\ Mumford.
  The statement here is due to \cite{Vidal} Proposition 2.4.2.1. 
(Cf.\ \cite{KKN1.5} Remark after {\sc Assumption} 8.1.)
  Another reference is \cite{Saito} Theorem 1.8. 
  Both papers based on the method of \cite{Yoshioka}.
(Actually, \cite{Yoshioka} and \cite{Saito} treat the case of log smooth fs log schemes over a discrete valuation ring, but the proof is in the same way. 
\cite{Saito} treats the non-vertical case also.)
  See \ref{TS} for a variant of Proposition \ref{desing}. 

  Finally, we discuss log de Rham cohomology and log Hodge theory. 

\begin{sbprop}
\label{prsmdR}
  Let $k$ be a field of characteristic zero. 
  Let $f \colon X \to S$ be a projective, log smooth and vertical morphism of fs log schemes with $S$ being log smooth over $k$.
  Let $q \in \bZ$. 
  Then we have the following. 

$(1)$ 
  $H^q_{\mathrm{dR}}(X/S):=R^q f_{\mathrm{k\acute{e}t}*}\omega_{X/S}^{\cdot,\mathrm{k\acute{e}t}}$ 
is a vector bundle endowed with a natural quasi-nilpotent integrable connection with log poles, 
and 
the Hodge filters 
$R^q f_{\mathrm{k\acute{e}t}*}\omega_{X/S}^{\cdot \ge p,\mathrm{k\acute{e}t}}$ are subbundles of 
$H^q_{\mathrm{dR}}(X/S)$ for all $p$.

$(2)$
When $k=\bC$, we have a natural log Hodge structure on $S_{\mathrm{k\acute{e}t}}$ of weight 
$q$ which is underlain by $H^q_{\mathrm{dR}}(X/S)$ with the Hodge filter. 
\end{sbprop}

\begin{pf}
  We may assume $k=\bC$, and (1) is deduced from (2). 
  (2) is by the main theorem of \cite{KMN} Theorem 8.1. %
\end{pf}

\begin{sblem}
\label{ketdR}
  Let $f\colon X \to S$ be a proper, log smooth and vertical morphism of fs log analytic spaces with $S$ being 
ideally log smooth over $\bC$ ({\rm \cite{IllusieKatoNakayama}} {\sc Definition (1.5)}). 
  Assume that for any $x$, the cokernel of $(M_S/\cO^{\times}_S)^{\mathrm{gp}}_{f(x)} \to (M_X/\cO^{\times}_X)^{\mathrm{gp}}_x$ is torsion-free. 
  Assume also that either $S$ is log smooth or $f$ is exact.
  Then we have a canonical isomorphism 
$$R^qf_*\omega^{\cdot,\mathrm{k\acute et}}_{X/S} = \varepsilon^*R^qf_*\omega^{\cdot}_{X/S}$$
for any $q\in \bZ$. 
  Here $\varepsilon$ is the forgetting-log morphism, i.e., the projection from the k\'et site to the usual site. 
\end{sblem}

\begin{pf}
  By \cite{IllusieKatoNakayama} {\sc Theorems} (6.2) and (6.3), 
the local system $R^qf_{*}^{\log}\bC$ corresponds to $R^qf_*\omega^{\cdot,\mathrm{k\acute et}}_{X/S}$ by the k\'et log Riemann--Hilbert correspondence, and 
it does to $R^qf_*\omega^{\cdot}_{X/S}$ by the non-k\'et log Riemann--Hilbert correspondence, 
respectively.
  Hence the desired isomorphism follows from the compatibility of the both Riemann--Hilbert correspondences (\cite{IllusieKatoNakayama} 
{\sc Theorem} (4.4)).
\end{pf}

\begin{sblem}
\label{dRblowup}
  Let the notation and the assumption be as in the previous lemma. 
  Let $X' \to X$ be a log blow-up and $f'\colon X' \to X \to S$ the composite. 
  Then the canonical homomorphism 
$$R^qf_*\omega^{\cdot}_{X/S}\to R^qf'_*\omega^{\cdot}_{X'/S}$$
is an isomorphism.
\end{sblem}

\begin{pf}
  By \cite{IllusieKatoNakayama} {\sc Theorem} (6.3), this homomorphism corresponds by the log Riemann--Hilbert correspondence 
to the homomorphism 
$R^qf^{\log}_*\bC\to R^qf^{\prime\log}_*\bC$
of local systems, which is an isomorphism by \cite{KajiwaraNakayama} {\sc Proposition} 5.3 (1).  
\end{pf}

\begin{sbprop}
\label{prsmdR2}
  Let $k$ be a field of characteristic zero. 
  Let $f \colon X \to s$ be a projective, log smooth and vertical morphism of fs log schemes with $s$ being 
the standard log point associated to $k$. 
  Let $q\in \bZ$. 
  Then we have the following. 

  $(1)$ $H^q_{\mathrm{dR}}(X/s):=R^q f_{\mathrm{k\acute{e}t}*}\omega_{X/s}^{\cdot,\mathrm{k\acute{e}t}}$ 
is a vector bundle with a natural quasi-nilpotent integrable connection with log poles. 

  $(2)$ When $k=\bC$, $H^q_{\mathrm{dR}}(X/s)$ 
carries a natural log Hodge structure on $s_{\mathrm{k\acute{e}t}}$ of weight $q$.
\end{sbprop}

\begin{pf}
  We may assume $k=\bC$, and (1) is deduced from (2). 
  We prove (2).  
  In \cite{FujisawaNakayama2}, the non-k\'et version of the case of (2) 
where $f$ is strict semistable is proved with the Hodge filter $R^q f_{*}\omega_{X/s}^{\cdot \ge p}$.
  We reduce (2) to this result as follows. 
  To prove (2), we slightly generalized the statement to the case where $s$ is the spectrum of a log Artin ring $\bC[\bN]/(x^n)$ for some 
$n \geq 1$, where $x$ is the generator of log. 
  We may assume that $f$ satisfies the assumptions in Lemma \ref{ketdR} by k\'et localization of the base $s$. 
  By a variant of Proposition \ref{desing}, we may assume further that 
there exists a log blow-up $X' \to X$ such that the special fiber of $X' \to s$ is strict semistable. 
  By Lemma \ref{ketdR}, we see that it is enough to show the non-k\'et version of (2).
  By the argument in \cite{IllusieKatoNakayama2} and the strict semistable case in \cite{FujisawaNakayama2}, $R^q f'_*\omega_{X'/s}^{\cdot}$ with the Hodge filters gives a log Hodge structure.  
  The non-k\'et version of (2) is reduced to this by Lemma \ref{dRblowup} 
and the induced Hodge filtration on 
$R^q f_{*}\omega_{X/s}^{\cdot}$ from $R^q f'_*\omega_{X'/s}^{\cdot}$ 
does not depend on the choice of $X'$. 
\end{pf}

\begin{sbprop}
\label{dRbasechange}
  Let $f \colon X \to S$ be a projective, log smooth and vertical morphism of fs log schemes with $S$ being log smooth over $\bC$.
  Let $s \to S$ be a standard log point associated to $\bC$ over $S$. 
  Let $f_s : X_s \to s$ be the base-changed morphism. 
  Let $q \in \bZ$. 
Then the pullback of the log Hodge structure $H^q_{\mathrm{dR}}(X/S)$ is naturally isomorphic to 
the log Hodge structure $H^q_{\mathrm{dR}}(X_s/s)$. 
\end{sbprop}

\begin{pf}
  Since there is a natural base change map, it is enough to show that the local system can be base-changed, 
which is by Proposition \ref{prsmbcB}.
\end{pf}

\subsection{Fans in log geometry}\label{fan}
Let $(\fs)$ be the category of fs log schemes which have charts Zariski locally. From now on, in the rest of this paper, an fs log scheme means an object of this $(\fs)$.

\medskip

We review the formulation of fans in \cite{KatoToric} as  unions of Spec of monoids.  This is a variant of the theory of polyhedral cone decompositions in \cite{KKMS}, \cite{Oda}. 

We add new things \ref{case2} and \ref{asfan2} which were not discussed in \cite{KatoToric}.

\begin{sbpara}
For a monoid $P$, an {\it ideal} of $P$ means a subset $I$ of $P$ 
such that $ab\in I$ for any $a\in P$ and $b\in I$. A {\it prime ideal} of $P$ means an ideal $\frak p$ of $P$ such that the complement $P\smallsetminus \frak p$ is a submonoid of $P$. We denote the set of all prime ideals of $P$ by $\Spec(P)$. 

\end{sbpara}

\begin{sbpara} For a monoid $P$ and for a submonoid $S$ of $P$, we have the monoid $S^{-1}P=\{s^{-1}a \;|\; a\in P, s \in S\}$ obtained from $P$ by inverting elements of $S$. Here $s_1^{-1}a_1=s_2^{-1}a_2$ if and only if there is an $s_3\in S$ such that $s_3s_2a_1=s_3s_1a_2$. 

In the case where $S=\{f^n\;|\; n\geq 0\}$ for $f\in P$, $S^{-1}P$ is denoted also by $P_f$.  
\end{sbpara}

\begin{sbpara} By a {\it monoidal space}, we mean a topological space $T$ endowed with a sheaf of monoids $\cP$ such that $(\cP_t)^\times=\{1\}$ for any $t\in T$. Here $\cP_t$ denotes the stalk of $\cP$ at $t$ and $(-)^\times$ means the subgroup consisting of all invertible elements. 
\end{sbpara}

\begin{sbpara} For a monoid $P$, $\Spec(P)$ is regarded as a monoidal space in the following way. 

We endow $\Spec(P)$ with the  topology for which the sets $D(f) = \{\frak p\in \Spec(P)\;|\; f\notin \frak p\}$ with $f\in P$ form a basis of open sets.

The sheaf $\cP$ of monoids on $\Spec(P)$ is characterized by the property that for $f\in P$, $\cP(D(f))= P_f/P_f^\times$.

The stalk of $\cP$ at $\frak p\in \Spec(P)$ is identified with $P_{\frak p}/(P_{\frak p})^\times$, where $P_{\frak p}= (P\smallsetminus \frak p)^{-1}P$.

\end{sbpara}

\begin{sbpara} For a monoidal space $\Sig$ with the structure sheaf $\cP$ of monoids and for a monoid $P$, the natural map
$\mathrm{Mor}(\Sig, \Spec(P)) \to \Hom(P, \cP(\Sig))$ is bijective. 

\end{sbpara}

\begin{sbpara}\label{deffan} A monoidal space  is called a {\it fan} if it has an open covering $(U_{\lambda})_{\lambda}$ such that each $U_{\lambda}$ is isomorphic, as a monoidal space, to $\Spec(P_{\lambda})$ for some monoid $P_{\lambda}$.

A fan which is isomorphic to $\Spec(P)$ for some monoid $P$ is called an {\it affine fan}. The functor $P\mapsto \Spec(P)$ is an anti-equivalence from the category of monoids $P$ such that $P^\times =\{1\}$ to the category of affine fans. The converse functor is given by $\Sig\mapsto \cP(\Sig)$, where $\cP$ is the structure sheaf of $\Sig$.
\end{sbpara}

\begin{sbpara}\label{XtoSig0}
For a fan $\Sig$, let 
$$[\Sig]:(\fs)\to (\mathrm{Sets})$$ 
be the contravariant functor which sends $X\in (\fs)$ to the set of all morphisms $(X, M_X/\cO_X^\times)\to \Sig$ of monoidal spaces. 

If $\Sig= \Spec(P)$, we have $[\Sig](X)= \Hom(P, \Gamma(X, M_X/\cO_X^\times))$. 

\end{sbpara}

\begin{sblem}\label{eqSig} The functor $\Sig\mapsto [\Sig]$ from the category of fans to the category of contravariant functors $(\fs)\to (\mathrm{Sets})$ is fully faithful. 

\end{sblem}

\begin{pf} This is reduced to the statement that the contravariant functor $P\mapsto [\Spec(P)]$ from the category of monoids $P$ such that $P^\times=\{1\}$ to the category of contravariant functors $(\fs)\to (\mathrm{Sets})$ is fully faithful. For monoids $P$ and $Q$ such that $P^\times=\{1\}$ and $Q^\times=\{1\}$ and for $X=\Spec(\Z[Q])$, we have $[\Spec(P)](X)=\Hom(P, \Gamma(X, M_X/\cO_X^\times))= \Hom(P,Q)$. From this, we obtain easily that the map  $\Hom(P, Q) \to \mathrm{Mor}([\Spec(Q)], [\Spec(P)])$ is bijective.  
\end{pf}

\begin{sbpara}\label{XtoSig} According to Lemma \ref{eqSig}, we will often identify a fan $\Sig$ with the functor $[\Sig]$. 

For an fs log scheme $X$ and for a fan $\Sig$, we will regard a morphism $(X, M_X/\cO_X^\times)\to \Sig$  of monoidal spaces as a morphism $X\to [\Sig]$ from the functor $X$ on $(\fs)$ represented by $X$ to the functor $[\Sig]$. We will call a morphism $X\to [\Sig]$ also a morphism $X\to \Sig$.  
\end{sbpara}

  The next lemma is easy to see. 

\begin{sblem}\label{strict} For an fs log scheme $X$, a fan $\Sig$, and a morphism $X\to \Sig$, the following conditions {\rm (i)} and {\rm (ii)} are equivalent.

{\rm (i)} The corresponding morphism $(X, M_X/\cO_X^\times) \to \Sig$ of monoidal spaces is strict. Here we say that a morphism $f: (T, \cP)\to (T', \cP')$ of monoidal spaces is {\rm strict} if $f^{-1}(\cP')\to \cP$ is an isomorphism. 

{\rm (ii)}  Locally on $X$, there is an open set $\Spec(P)$ of $\Sig$ with $P$ a monoid such that $X\to \Sig$ factors as  $X\to \Spec(\Z[P]) \to \Spec(P) \subset \Sig$, where $\Spec(\Z[P])$ is endowed with the standard log structure and the homomorphism $P\to M_X$ corresponding to the first arrow is a chart of $X$ (that is, the first morphism is strict, where we say a morphism of log schemes  $X\to Y$ is {\rm strict} if the log structure of $X$ coincides with the inverse image of the log structure of $Y$). 
\end{sblem} 

\begin{sbpara}
We will say $X\to \Sig$ is {\it strict} if the equivalent conditions in Lemma \ref{strict} are satisfied.
\end{sbpara}

\begin{sbpara}\label{clfan} 
  Polyhedral cone decompositions  which appear in the toric geometry 
(\cite{KKMS}, \cite{Oda}) are related to the above notion of fan (\ref{deffan}) as follows.

Let $N$ be a free $\Z$-module of finite rank,
and let $N_\R := \R \otimes_{\Z} N$.
A {\it rational polyhedral cone} in $N_{\R}$ is a subset
of the form 
$$ \sigma = \big\{ \sum_{i=1}^r x_i N_i \ \big| \ x_i \in \R_{\geq 0} \big\} $$
for some $N_1,\dots,N_r \in N$.
A rational polyhedral cone $\sigma$ is called {\it strongly convex}
if it does not contain a line, i.e.\ $\sigma \cap (-\sigma) = \{ 0 \}$.
A subset $\tau \subset \sigma$ is called a {\it face} of $\sigma$
if there exists an element $h\in \Hom_{\R}(N_\R, \R)$
such that $\sigma \subset \{ x \in N_{\R} \mid h(x) \geq 0 \}$
and $\tau = \sigma \cap \{ x \in N_{\R} \mid h(x) = 0 \}$. A face of $\sig$ is also a rational polyhedral cone.
\medskip

A rational polyhedral cone decomposition in $N_\R$ (or a  rational fan in $N_\R$) is a non-empty 
set $\Sig$
of strongly convex rational polyhedral cones in $N_\R$
satisfying the following two conditions:
(i) If $\sig \in \Sig$ and $\tau$ is a face of $\sig$, then $\tau \in \Sig$;
(ii) If $\sig, \tau\in \Sig$, the intersection $\sig \cap \tau$ is a face of $\sig$.

\medskip
We regard a rational fan $\Sig$ in $N_{\R}$ as a fan in the sense of \ref{deffan} as follows. 

We endow $\Sig$ with the topology for which the sets $\mathrm{face}(\sig)$ of all faces of $\sig$ for  $\sig\in \Sig$ form a basis of open sets. 

We endow $\Sig$ with the sheaf $\cP$ of monoids characterized by $\cP(\mathrm{face}(\sig))= P_{\sig}/(P_{\sig})^\times$, where  
$$P_{\sig}= \{h\in \Hom(N, \Z)\;|\; h(x) \geq 0\;\text{for all}\;x\in \sig\}.$$

The open set $\mathrm{face}(\sig)$ of $\Sig$ is identified with $\Spec(P_{\sig})$. 
 \end{sbpara}
 
 \begin{sbpara}\label{toric}
   For a rational fan $\Sig$ in $N_\R$, we have the toric variety $\mathrm{Toric}_{\Sig}= \bigcup_{\sig\in \Sig} \Spec(\Z[P_{\sig}])$ over $\Z$ 
 corresponding to $\Sig$ with the standard log structure, on which the torus $N \otimes {\bf G}_m$ acts naturally. We have
 $$[\Sig]= \mathrm{Toric}_{\Sig}/(N \otimes {\bf G}_m)$$ as a sheaf on $(\fs)$, where  $\mathrm{Toric}_{\Sig}$ is identified with the sheaf on $(\fs)$ represented by it. 
\end{sbpara}

\begin{sbpara}\label{asfan}

For an fs log scheme $X$, in the following cases (i) and (ii),  we can define a fan $\Sig_X$ associated to $X$ and a strict morphism $X\to \Sig_X$ in a canonical way. 

\medskip

Case (i). $X$ is log regular (\cite{KatoToric}).

Case (ii). $X$ is vertical and log smooth over a standard log point. 

\medskip

The case (i) was considered in \cite{KatoToric}. The case (ii) is explained below. 
\end{sbpara}

\begin{sbpara}\label{case1} We review first the case (i). See \cite{KatoToric} for the definition of log regularity. As a set, $\Sig_X$ is the set of all points $x$ of $X$ such that the maximal ideal $m_x$ of $\cO_{X,x}$ is generated by the image of $M_{X,x}\smallsetminus \cO^\times_{X,x}$, where $M_{X,x}$ is the stalk at $x$ of the direct image of $M_X$ to the Zariski site. 
The topology of $\Sig_X$ is the restriction of the topology of $X$. The structural sheaf $\cP$ of $\Sig_X$ is defined as the inverse image of the sheaf 
$M_X/\cO_X^\times$ on $X$. The morphism $(X, M_X/\cO_X^\times)\to \Sig_X$ is defined as follows. As a map, it sends $x\in X$ to the point of $X$ corresponding to the prime ideal of $\cO_{X,x}$ generated by the image of $M_{X,x}\smallsetminus \cO^\times_{X,x}$. If $x\in X$ and if $y\in X$ 
is the image of $x$ in $\Sig_X$, there is a chart $P\to M_U$ for some open neighborhood $U$ of $x$ in $X$ such that $P\to (M_X/\cO_X^\times)_x$ is an isomorphism, 
and via the composite homomorphism $P\to (M_X/\cO_X^\times)_x\to  (M_X/\cO_X^\times)_y$, $\Spec(P)$ is identified with an open neighborhood of $y$ in $\Sig_X$. The chart defines a morphism $(U, M_U/\cO_U^\times) \to \Spec(P)$ and hence a morphism $(U, M_U/\cO_U^\times) \to \Sig_X$ and these local definitions glue to a global definition of $(X, M_X/\cO_X^\times) \to \Sig_X$.

\end{sbpara}

\begin{sbpara}\label{case2} We consider the case (ii). As a set, $\Sig_X$ is the disjoint union $\Sig_X'\coprod \{\eta\}$ of the set $\Sig_X'$ of all points $x$ of $X$ such that the maximal ideal $m_x$ of $\cO_{X,x}$
 is generated by the image of $M_{X,x}\smallsetminus \cO^\times_{X,x}$ and the one-point-set $\{\eta\}$. The topology on $\Sig_X$ is as follows. First define the topology
  of $\Sig'_X$ to be the restriction of the topology of $X$. A closed subset of $\Sig_X$ is either a closed subset of $\Sig'_X$ or $\Sig_X$. 
  The structure sheaf $\cP$ of monoids on $\Sig_X$ is defined as follows. First let the sheaf $\cP'$ on $\Sig'_X$ be the inverse image of $M_X/\cO_X^\times$. Let $\cP=i_*\cP'$, where 
$i:\Sig'_X \to X$ is the inclusion map. 

Then $\Sig_X$ is a fan. 
  This is reduced to the log regular case, since $X$ is locally isomorphic to the part $t=0$ of a log regular scheme $Y$, where $t$ is a section of log structure $M_Y$ of $Y$ such that the part of $Y$ where $t$ is invertible coincides with the part where $M_Y$ is 
trivial.

We define a map $X\to \Sig_X$ in the similar way to the case (i) described above.  
  The proof for the gluing also reduces to the case (i). 
  The resulting map in fact factors through $X\to \Sig'_X$. 
\end{sbpara}

\begin{sbpara}\label{asfan2} Outside the cases (i) and (ii) in \ref{asfan}, it seems difficult to develop a general theory of  fans canonically associated to fs log schemes (cf.\ \cite{ACMUW}). 
We give an example of an fs log scheme  $X$ having the following nice property (1) but such that for any fan $\Sig$,  there is no  strict morphism $(X, M_X/\cO^\times_X)\to \Sig$.

(1)  $X$ is locally isomorphic to a closed subscheme of a log regular scheme $Y$ defined by an ideal of $\cO_Y$ generated by the images of sections of the log structure  $M_Y$ of $Y$ under $M_Y\to \cO_Y$ endowed with the log structure induced by the log structure of $Y$. 
As a scheme, $X$ is a union of two ${\bf P}^1_k$ obtained by identifying  $0$ of each ${\bf P}^1$  with $\infty$ of the other ${\bf P}^1$.

Let $k$ be a field. Endow $\Spec(k[x_1,x_2,x_3, x_4])$ with the log structure associated to $\N^4\to k[x_1,x_2,x_3,x_4]\;;\; n\mapsto \prod_{i=1}^4 x_i^{n(i)}$. Let $Z= \Spec(k[x_1, x_2, x_3, x_4]/(x_1x_2, x_3, x_4))$ with the induced log structure, 
 and let $Z'$ be a copy of $Z$. (Hence as schemes, $Z$ and $Z'$ are isomorphic to $\Spec(k[x,y]/(xy))$.) Denote the copy of $x_i$  on $Z'$ by $x'_i$. Let $U$ be the part of $Z$ on which $x_1$ is invertible and let $V$ be the part of $Z$ on which $x_2$ is invertible. Let $U'$ and $V'$ be the copies of $U$ and $V$ in $Z'$, respectively. Let $X$ be the union of $Z$ and $Z'$ which we glue by identifying the open set $U \coprod V$ of $Z$ and the open set $U'\coprod V'$ of $Z'$, as follows. We identify $U$ and $U'$ 
by identifying $x'_1$ with $1/x_1$, $x'_2$ with $x_1^2x_2$, $x'_3$ with $x_3$, and $x'_4$ with $x_4$  in the log structure. (Hence $x_2$ is identified with $(x'_1)^2x_2'$ in the log structure.) We identify 
$V$ and $V'$ by identifying $x'_2$ with $1/x_2$, $x'_1$ with $x_1x_2^2$, $x'_3$ with $x_4$, and $x'_4$ with $x_3$  in the log structure. (Hence $x_1$ is identified with $x'_1(x'_2)^2$ in the log structure.) 

We show that there is no strict morphism $f: X\to \Sig$ to any fan $\Sig$. 

Assume $f$ exists. Let $p$ be the point of $Z$ at which all $x_i$ have value $0$, let $p'\in Z'$ be the copy of $p$, let $u$ be the generic point of $U$, and let $v$ be the generic point of $V$. Let $\cP$ be the structure sheaf of monoids of $\Sig$. Then $\cP_{f(p)}$ is identified with $(M_X/\cO_X^\times)_p\cong \N^4$ which is generated by $x_1,x_2,x_3,x_4$. $\Sig$ has an open neighborhood which is identified with $\Spec(\cP_{f(p)})$. Since $p$ belongs to the closure of $u$ in $X$, $f(u)$ belongs to $\Spec(\cP_{f(p)})$. We have a commutative diagram 
$$\begin{matrix}  \cP_{f(p)}  &\to & \cP_{f(u)} \\ \downarrow && \downarrow\\
(M_X/\cO_X^\times)_p &\to &(M_X/\cO_X^\times)_u\end{matrix}$$
in which vertical homomorphisms are isomorphisms,
and hence $f(u)$ is the prime ideal of $\cP_{f(p)}$ generated by $x_2, x_3, x_4$. The open neighborhood of $u$ in $\Sig$ which is identified with $\Spec(\cP_{f(u)})$ is regarded as an open set of $\Spec(\cP_{f(p)})$. In this identification, the prime ideal of $\cP_{f(p)}$ generated by $x_3$ is identified with the prime ideal of $\cP_{f(u)}$ generated by $x_3$. Similarly, $\Spec(\cP_{f(u)})$ is identified with an open set of $\Spec(\cP_{f(p')})$ and the prime ideal of $\cP_{f(u)}$ generated by $x_3$ is identified with the prime ideal of $\cP_{f(p')}$ generated by $x_3'$.

Similarly $\Spec(\cP_{f(v)})$ is identified with an open set of $\Spec(\cP_{f(p)})$ and also with an open set of $\Spec(\cP_{f(p')})$. The prime ideal of $\cP_{f(v)}$ generated by $x_4$ is identified with the prime ideal of $\cP_{f(p)}$ generated by $x_4$ and it is also identified with the prime ideal of $\cP_{f(p')}$ generated by $x'_3$. This shows that 
the prime ideal of $\cP_{f(p)}$ generated by $x_3$ is equal to the prime ideal generated by $x_4$. Contradiction. 

\end{sbpara}

\subsection{Subdivisions of fans and log modifications}
\label{logmodification}
\begin{sbpara}
We shall mainly consider fans $\Sig$ (\ref{deffan}) satisfying the following condition $(S_{\mathrm{fan}})$ (like in \cite{KatoToric}). 

$(S_{\mathrm{fan}})$\; There exists an open covering $(U_{\lambda})_{\lambda}$ such that for each $\lambda$, $U_{\lambda}\cong \Spec(P_{\lam})$ as a fan for some fs monoid $P_{\lam}$. 
\end{sbpara}

\begin{sbpara} Let $N$ be as in \ref{clfan}, let $\sig$ be a strictly convex rational polyhedral cone in $N_\R$, and let 
$\Sig$ be the rational fan $\mathrm{face}(\sig)$ in $N_\R$ consisting of all faces of $\sig$. Then a {\it finite subdivision} of $\Sig$ means a finite rational fan $\Sig'$ in $N_{\R}$ such that $\sig=\bigcup_{\tau\in \Sig'} \tau$.

\end{sbpara}

\begin{sblem}\label{fsd1}  Let $\Sig=(\Sig, \cP)$ and $\Sig'=(\Sig', \cP')$ be fans satisfying the condition $(S_{\mathrm{fan}})$ and let $f:\Sig'\to \Sig$ be a morphism of fans. Then the following conditions {\rm (i)} and {\rm (ii)} are equivalent.  

{\rm (i)} $f$ satisfies the following {\rm (i-1)}--{\rm (i-3)}. 

{\rm (i-1)}  For any $t\in \Sig$, the inverse image $f^{-1}(t)$ is finite. 

{\rm (i-2)} For any $t\in \Sig'$, $\cP^{\gp}_{f(t)}\to  (\cP')_t^{\gp}$ is surjective. 

{\rm (i-3)} The map $\mathrm{Mor}(\Spec(\N), \Sig') \to \mathrm{Mor}(\Spec(\N), \Sig)$ is bijective.

{\rm (ii)} There exists an open covering $(U_{\lam})_{\lam}$ of $\Sig$ such that for each $\lam$, there are a finitely generated free $\Z$-module $N_{\lam}$, a strongly convex rational polyhedral cone $\sig_{\lam}$  in $N_{\lam, \R}$, a finite subdivision $V_{\lam}$ of 
$\mathrm{face}(\sig_{\lam})$, 
and a commutative diagram of fans 
$$\begin{matrix}  U'_{\lam}  &\cong& V_{\lam}\\
\downarrow &&\downarrow \\
U_{\lam} &\cong & \mathrm{face}(\sig_{\lam}),\end{matrix}$$
where $U'_{\lam}$ denotes the inverse image of $U_{\lam}$ in $\Sig'$. 

\end{sblem}

\begin{pf}
The proof is straightforwards. 
\end{pf}

\begin{sbpara}\label{fsd2} Let $\Sig$ be a fan satisfying $(S_{\mathrm{fan}})$. A {\it finite subdivision} of $\Sig$ (called a proper subdivision of $\Sig$ in \cite{KatoToric})  is a fan $\Sig'$ satisfying $(S_{\mathrm{fan}})$ endowed with a morphism $\Sig'\to \Sig$ satisfying the equivalent conditions (i) and (ii) in Lemma \ref{fsd1}. 

\end{sbpara}

\begin{sblem} Let $\Sig$ be a fan satisfying the condition $(S_{\mathrm{fan}})$, let $X$ be an fs log scheme, 
let $X\to \Sig$ be a morphism ($\ref{XtoSig}$), and let $\Sig'$ be a  finite subdivision of $\Sig$. Then the functor $X\times_{\Sig} \Sig': (\fs)\to (\mathrm{Sets})$ is represented by an fs log scheme $X'$ which is proper and log \'etale over $X$. Here $X\times_{\Sig} \Sig'$ denotes the fiber product of the functors $X=\Mor(\;, X)$ and $\Sig'=[\Sig']$ ($\ref{XtoSig0}$) on $(\fs)$ over the functor $\Sig= [\Sig]$ on $(\fs)$ (it does not mean the set theoretic fiber product of $X$ and $\Sig'$ over $\Sig$).

\end{sblem}

\begin{pf} We are reduced to the case $\Sig= \mathrm{face}(\sig)$ for a strongly convex rational polyhedral cone $\sig$ and $\Sig'$ is a finite subdivision of $\Sig$. Locally on $X$, $X\to \Sig$ is the composition 
$X\to \Spec(\Z[P_{\sig}])\to \Sig$. Hence we are reduced to the case 
$X=\Spec(\Z[P_{\sig}])$. Then
 $X \times_{\Sig} \Sig'$ is represented by the toric variety $\bigcup_{\tau\in \Sig'} \Spec(\Z[P_{\tau}])$ over $\Z$ 
associated to $\Sig'$, which is proper and log \'etale over $X$. 
\end{pf}

\begin{sbpara}
\label{logmod}
We call a morphism $X\to Y$ of fs log schemes a {\it log modification} if locally on $Y$, there are a fan $\Sig$ satisfying $(S_{\mathrm{fan}})$, a morphism $Y\to \Sig$, and a  
 finite subdivision $\Sig'$ of $\Sig$ such that $X$ represents $Y \times_{\Sig} \Sig'$. 
  
\end{sbpara}

Log modifications were studied in \cite{KatoUsui} for fs log analytic spaces over $\C$. 

The following lemma is easy to prove.
\begin{sblem}\label{lemlm}

 $(1)$ A log modification is proper and log \'etale. 

 $(2)$ If $X\to Y$ is a log modification, the induced morphism of functors $\mathrm{Mor}(\;,X) \to \mathrm{Mor}(\;,Y)$ on $(\fs)$ is injective.

 $(3)$ If $X_i \to Y$ ($i=1,2$) are log modifications, $X_1 \times_Y X_2 \to Y$ is a log modification. Here $X_1\times_Y X_2$ denotes the fiber product in the category of fs log schemes. 
 
 $(4)$ If $X\to Y$ and $Y\to Z$ are log modifications, the composition $X\to Z$ is a log modification.
 \end{sblem}

\begin{sbprop}\label{invet}
  Let $f: X \to Y$ be a log modification of fs log schemes.

$(1)$ Let $F$ be a torsion sheaf of abelian groups on 
$Y_{\mathrm{l\acute{e}t}}$. 
  Then the natural homomorphism $F \to Rf_{\mathrm{l\acute{e}t}*}f^*_{\mathrm{l\acute{e}t}}F$ is an isomorphism. 

$(2)$ Let $\ell$ be a prime number which is invertible on $Y$. 
  Then the natural homomorphism 
$\bQ_{\ell} \to Rf_{\mathrm{l\acute{e}t}*}f^*_{\mathrm{l\acute{e}t}}\bQ_{\ell}$ is an isomorphism. 
\end{sbprop}

\begin{pf}
  (2) is reduced to (1). 
  (1) is a slight generalization of Theorem 5.5 (2) of \cite{Nakayama3}, and the proof is similar, which is reduced easily to Lemma \ref{lemlm} (2). 
\end{pf}

\begin{sbpara}\label{moreg} (1) Let $\Sig$ be a fan with the structure sheaf $\cP$  of monoids. We say $\Sig$ is {\it free} if for any $t\in \Sig$, the stalk $\cP_t$ is isomorphic to $\N^{r(t)}$ for some $r(t)\geq 0$. 

(2) Let $X$ be an fs log scheme. We say $M_X/\cO_X^\times$ is {\it free} if for any $x\in X$, $(M/\cO_X^\times)_x\cong \N^{r(x)}$ for some $r(x)\geq 0$. 
\end{sbpara}

\begin{sbprop}\label{Propdesing} Let $\Sig$ be a finite fan satisfying the condition $(S_{\mathrm{fan}})$. Then there is a finite subdivision $\Sig'\to \Sig$ which is free ($\ref{moreg}$ $(1)$). 
\end{sbprop}

This is already explained in \cite{KatoToric}.

\begin{sblem}\label{extsd}
Let $\Sig$  be a finite fan satisfying the condition $(S_{\mathrm{fan}})$ with the structural sheaf $\cP$, let $t\in \Sig$, and let $P$ be an fs submonoid of $\cP^{\gp}_t$ containing $\cP_t$. Then  there is a finite subdivision $\Sig'$ of $\Sig$ such that there is an open immersion $\Spec(P)\to \Sig'$ over $\Sig$. 
\end{sblem}

\begin{pf} 
  Regard $\Sigma$ as a conical polyhedral complex with an integral structure 
(\cite{KKMS} Chapter II, \S1, Definitions 5 and 6, pp.69--70).
  Let $\sig$ be its cell corresponding to $\cP_t$ and $\tau \subset \sigma$ be the subcone corresponding to $P$.
  Take a rational homomorphism $f\colon \sig \to \bR_{\ge0}$ such that 
$f^{-1}(\{0\})$ is trivial, 
where $\bR_{\ge0}$ is the monoid of the nonnegative real numbers with addition. 
  Let $f_0\colon S:=\underset {\sigma' \in \Sigma} \bigcup \mathrm{Sk}^{1}(\sigma') \cup \mathrm{Sk}^{1}(\tau) \to \bR$ be the 
zero extension of the restriction of $f$ to $\mathrm{Sk}^{1}(\tau)$, 
that is, for any $s \in S$, $f_0(s)=f(s)$ if $s \in \mathrm{Sk}^{1}(\tau)$ and $f_0(s)=0$ otherwise.
Here $\mathrm{Sk}^1$ means the 1-skeleton 
(\cite{KKMS} Chapter I, \S2, p.29). 
  Let $f_1\colon |\Sigma| \to \bR_{\ge0}$ be the convex interpolation of $f_0$ (\cite{KKMS}, Chapter I, \S2, p.29 and Chapter II, \S2, p.92), where $|\Sigma|$ is the support of $\Sigma$. 
  Then, $f_1$ coincides with $f$ on $\tau$, and the coarsest subdivision 
of the conical polyhedral complex $\Sigma$ on any cell of which $f_1$ is linear owes $\tau$ as a cell.
  Hence the corresponding finite subdivision $\Sig'$ of the fan $\Sigma$ satisfies the desired property. 
\end{pf}

\begin{sbprop}\label{cfan}
 Let $X$ be a quasi-compact fs log scheme, let $\Sig$ be a finite fan satisfying the condition $(S_{\mathrm{fan}})$ with the structure sheaf $\cP$, and let $f:X\to \Sig$ be a morphism ($\ref{XtoSig}$) such that for any $x\in X$, the map $\cP_{f(x)}\to (M_X/\cO_X^\times)_x$ is surjective. Then for a sufficiently fine finite
subdivision $\Sig'$ of $\Sig$, $X\times_{\Sig} \Sig'\to \Sig'$
is strict.
\end{sbprop}

\begin{pf}
  First notice that the problem is local on $X$ as the category of finite subdivisions of $\Sigma$ is directed. 
Let $x\in X$, and let
$P$ be the fs submonoid of $(\cP_{f(x)})^{\gp}$ consisting of all
elements whose images in $(M^{\gp}_X/\cO_X^\times)_x$ are contained in
$(M_X/\cO_X^\times)_x$. Then $P/P^\times\to (M_X/\cO^\times_X)_x$ is
an isomorphism. 
  Since $X$ is quasi-compact and the problem is local on $X$, replacing $X$ by an open neighborhood of $x$, we may assume that $X \to \Sig$ factors 
as $X\to \Spec(P) \to \Spec(\cP_{f(x)}) \to \Sigma$ and the first arrow is strict. Let $\Sig'$ be a finite subdivision of $\Sig$ such  that there is an open immersion $\Spec(P)\to 
\Sig'$ over $\Sig$ (Lemma \ref{extsd}). Then the morphism $X= X \times_{\Sig} \Sig' \to \Sig'$ is strict because it is the composition of strict morphisms $X \to \Spec(P) \to \Sig'$. 
\end{pf}

\begin{sbrem} This Proposition \ref{cfan} will be used later in Proposition \ref{iddia1} to make the diagonal of a vertical log smooth fs log scheme over a
standard log point a regular immersion, by log modification. 
\end{sbrem}

\begin{sbpara}
\label{TS} 
  From the next section, we will use the following corollary of Proposition \ref{desing}.

Let $X$ be a projective vertical log smooth fs log scheme over a standard log point $s$. %
  Then, for some morphism of standard log points $s' \to s$ whose underlying extension of the fields is 
an isomorphism, %
  we have a semistable fs log scheme $X'$ over $s'$ which is a log blow-up of $X \times_ss'$. %
\end{sbpara}

\subsection{Grothendieck groups of vector bundles and log geometry}
\label{Grothendieckgp}
\begin{sbpara}
  Recall the following theory in \cite{SGA6} till \ref{Kpush}. 

For a scheme $X$, let $K(X)$ be the Grothendieck group of
the category of locally free $\cO_X$-modules on $X$ of finite rank. It is a commutative ring in which the multiplication corresponds to tensor products.

The $K$-group $K(X)$ has a decreasing filtration $(F^r K(X))_{r \in \Z}$
called the {\it $\gamma$-filtration}
(for details, see \cite{SGA6}, \cite{FultonLang}, Chapter III, V).
It satisfies
$F^0 K(X) = K(X)$ and
$F^r K(X) \cdot F^s K(X) \subset F^{r+s} K(X)$.
We define
$$ \gr^r K(X) := F^r K(X) / F^{r+1} K(X). $$
\end{sbpara}

\begin{sbpara}
\label{Kpush}
For a morphism $X\to Y$ of schemes, the pullback homomorphism $K(Y)\to K(X)$ is defined and it respects the $\gamma$-filtration. 

On the other hand, for a morphism $f: X\to Y$ of schemes which is projective and locally of complete intersection (cf.\ \cite{SGA6}, Expos\'e VIII, D\'efinition 1.1),  the pushforward homomorphism $K(X)\to K(Y)$ is defined 
(cf.\ \cite{SGA6} Expos\'e IV, 2.12). 
It sends $F^iK(X)_{\Q}$ to $F^{i-d}K(Y)_{\Q}$. 
  Here $d$ is the relative dimension of $f$ which is a locally constant function on $X$ characterized as follows. Locally on $X$,  $f$ is a composition $X\overset{i}\to Z \overset{g}\to Y$, where $i$ is a regular immersion and $g$ is smooth.  The relative dimension of $f$ is $d_1-d_2$, where $d_1$ is the relative dimension of $g$ and $d_2$ is the  codimension of $i$. 
\end{sbpara}

\begin{sbpara} If $X$ and $Y$ are projective smooth schemes over a field $k$, any morphism $X\to Y$ over $k$ is projective and locally of complete intersection and hence the pushforward homomorphism $K(X) \to K(Y)$ is defined. 
However, in log geometry, we have no such a nice property if we replace the smoothness by  log smoothness.

We give some preliminaries to treat log smooth situations which we encounter in later sections. 
\end{sbpara}

\begin{sbprop}\label{flat} Let $S$ be an fs log scheme of log rank $\leq 1$ (this means that for any $s\in S$, $(M_S/\cO_S^\times)_s$ is isomorphic to either $\N$ or $\{1\}$). 
Let $f:X\to S$ be a log smooth morphism. Then the underlying morphism of schemes of $f$ is flat. 
\end{sbprop}

\begin{pf}
\cite{KatoLog}, {\sc Corollary} (4.4), {\sc Corollary} (4.5).
\end{pf}

\begin{sbprop}
\label{PropLci}
Let $S$ be an fs log scheme of log rank $\leq 1$,
and let $f \colon X \to Y$ be a morphism of fs log schemes over $S$.
Assume that $X,Y$ are log smooth over $S$,
and assume that 
 $M_X/\cO_X^\times$ and $M_Y/\cO_Y^\times$ are free ($\ref{moreg}$). 
Then the underlying morphism of schemes of $f$
is locally of complete intersection.
\end{sbprop}

\begin{pf} Working \'etale locally on $X$ and on $Y$, we may assume that $f$ is the base change of $f': X'\to Y'$ over $S'=\Spec(\Z[\bN])$ by a strict morphism $S\to S'$, where $S'$ is endowed with log by $\bN$ and $X'$ and $Y'$ are log smooth over $S'$. By the assumption on the log of $X$ and $Y$, we may assume that 
$M/\cO^\times$ of $X'$ and that of $Y'$ are also free (\ref{moreg})
and hence $X'$ and $Y'$ are smooth over $\Z$ as schemes. Hence $f'$ is locally of complete intersection.
  Since $X'$ and $Y'$ are flat over $S'$, $f$ is also locally of complete intersection.  
  Here we used the fact that any base change of a morphism 
$f'\colon X' \to Y'$ of locally complete intersection of schemes which are flat over a scheme is locally of complete intersection. 
  A proof of this fact is as follows. 
  Locally, $f'$ is the composition of a regular immersion followed by a smooth morphism, and hence we may assume that $f'$ is a regular immersion. But for a closed immersion defined by an ideal $I$ being a regular immersion is equivalent to the condition that $I/I^2$ is locally free and $I^n/I^{n+1}=\mathrm{Sym}^n(I/I^2)$ for any $n$. The last property is stable under any base change. %
\end{pf}

\begin{sbpara} For an fs log scheme $X$, we define
$$K_{\lim}(X):=\varinjlim_{X'} K(X'),$$
where $X'$ ranges over all log modifications (\ref{logmod}) of $X$.
\end{sbpara}

\begin{sblem} \label{tildelem}
 Let $X$ be a quasi-compact fs log scheme, let $\Sig$ be a finite fan satisfying the condition $(S_{\mathrm{fan}})$ with the structure sheaf $\cP$, and let $f:X\to \Sig$ be a morphism ($\ref{XtoSig}$) such that for any $x\in X$, the map $\cP_{f(x)}\to (M_X/\cO_X^\times)_x$ is surjective. Then 
 we have an isomorphism
 $$\varinjlim_{\Sig'} K(X\times_{\Sig} \Sig')\overset{\cong}\to K_{\lim}(X),$$
where $\Sig'$ ranges over all finite subdivisions of $\Sig$. 

\end{sblem}

\begin{pf} Let $X'\to X$ be a log modification. Then the composition $f': X' \to X \to \Sig$  satisfies the condition that 
 $\cP_{f'(x)}\to (M_{X'}/\cO_{X'}^\times)_x$ is surjective for any $x\in X'$. Hence by Proposition \ref{cfan}, there is a finite subdivision $\Sig'$ of $\Sig$ such that the morphisms $X\times_{\Sig} \Sig' \to \Sig$ and $X' \times_{\Sig} \Sig'\to \Sig'$ are strict. This shows that the log modification $X' \times_{\Sig} \Sig'\to X\times_{\Sig} \Sig'$ is strict and hence $X' \times_{\Sig} \Sig'\overset{\cong}\to X\times_{\Sig} \Sig'$. 
\end{pf}

\begin{sbpara}\label{Chern} Let $s$ be a geometric standard log point (\ref{slp}), 
and let $X$ be an fs log scheme over $s$. Let $\ell$ be a prime number which is different from the characteristic of $s$ and let  $H^m(X)_{\ell}:=R^mf_*\Q_{\ell}$, where $f$ is the morphism $X\to s$ and $R^mf_*$ is the $m$-th higher direct image for the  log \'etale topology (\ref{smoothQ_l}). We will identify $H^m(X)_{\ell}$  with its stalk. 

We have a Chern class map
$\gr^i K(X)_\Q \to H^{2i}_{\mathrm{\acute{e}t}}(X, \Q_{\ell})(i)$ to the classical \'etale cohomology, which coincides with the Chern character map. By composing this with the canonical map $H^{2i}_{\mathrm{\acute{e}t}}(X, \Q_{\ell})(i)
\to H^{2i}(X)_{\ell}(i)$ and by going to the inductive limit for log modifications using the invariance Proposition \ref{invet} for the log \'etale cohomology, we obtain the Chern class map
$$ \gr^i K_{\lim}(X)_\Q\to H^{2i}(X)_{\ell}(i).$$
\end{sbpara}

\begin{sbprop}\label{RR}
Let $X$ (resp.\  $Y$) be projective and vertical log smooth fs log scheme over a geometric standard log point $s$ such that ($\ref{slp}$) $M/\cO^\times$ of $X$ and that of $Y$ are free ($\ref{moreg}$). 
Let $f:X\to Y$ be a morphism over $s$ of relative dimension $d$. ($d$ can be $< 0$. Cf.\ $\ref{Kpush}$.) Let $\ell$ be a prime number which is different from the characteristic of $s$. Then for any $i\in \Z$, the following diagram is commutative. 
$$\begin{matrix} \gr^{i+d} K (X)_\Q & \to & H^{2(i+d)}(X)_{\ell}(i+d)\\
\downarrow  &&\downarrow\\
\gr^iK(Y)_\Q  &\to & H^{2i}(Y)_{\ell}(i).
\end{matrix}$$
Here the left vertical arrow is defined by Proposition $\ref{PropLci}$ and $\ref{Kpush}$ and the right vertical arrow is the pushforward map (the dual of $H^{2j}(Y)_{\ell}(j)\to H^{2j}(X)_{\ell}(j)$ 
for Poincar\'e duality (Corollary $\ref{Q_lPoin}$), where $j=\dim(Y)-i$).
\end{sbprop}

\noindent {\it Remark.}\ 
  In the above, $d$ (resp.\   $\dim(Y)$) is considered as a locally constant function on $X$ (resp.\  $Y$) (cf.\ \ref{Kpush}).
  In general, if $m$ is a locally constant function on $X$, $H^m(X)$ means $\bigoplus_i \; H^{m(i)}(X_i)$, where $X_i$ are connected components of $X$ and $m(i)$ is the value of $m$ on $X_i$. The meaning of $\gr^m K(X)_\Q$ is similar. 

\begin{pf} Let $X^{\circ}$ (resp.\  $Y^{\circ}$) be the underlying scheme of $X$ (resp.\  $Y$).  The morphism $f$ is the composition of two morphisms $X\to  {\bf P}^n \times Y \to Y$ in which the underlying morphism of schemes of the first arrow is a closed regular immersion and the second arrow is the projection.  It is sufficient to prove Proposition \ref{RR} for each of these two morphisms. The proof for the latter morphism is standard. We consider the first morphism. It is sufficient to
 prove the commutativity of the two squares in the diagram 
$$\begin{matrix} \gr^i K(X^{\circ})_\Q &\to & H^{2i}_{\mathrm{\acute{e}t}}(X^{\circ}, \Q_{\ell})(i)& \to & H^{2i}(X)_{\ell}(i)\\
\downarrow && \downarrow && \downarrow \\
 \gr^{i+c} K(Y^{\circ})_\Q &\to & H^{2i+2c}_{\mathrm{\acute{e}t}}(Y^{\circ}, \Q_{\ell})(i+c)& \to & H^{2i+2c}(Y)_{\ell}(i+c)\end{matrix}$$
 assuming that the morphism $X^{\circ}\to Y^{\circ}$ is a closed regular immersion of codimension $c$. 
  Here 
the central vertical arrow is the Gysin map which is defined as follows.  Let $\xi\in H^{2c}_{X^{\circ}}(Y^{\circ}, \Q_{\ell})(c)$ be the localized Chern class of the $\cO_Y$-module $\cO_X$ (\cite{Iversen}). 
By using the cup product 
$$H^i_{\mathrm{\acute{e}t}}(X^{\circ}, \Q_{\ell}) \times H^j_{X^{\circ}}(Y^{\circ}, \Q_{\ell})\to H^{i+j}_{X^{\circ}}(Y^{\circ}, \Q_{\ell}),$$
let the Gysin map be the product with $\xi$. 
(Cf. \cite{BFM} Section 5.4.)

The left square is commutative by Riemann--Roch theorem in 
Corollary 1 in Section 5.3 of \cite{BFM} (see also \cite{Fulton}). 
We prove that the right square is commutative. By \ref{TS}, we may assume that $X$ and $Y$ are semistable. 
Let $X'$ be $X^{\circ}$ with the inverse image of the log structure of $Y$. Hence $X\to Y$ factors as $X\to X'\to Y$. Consider the diagram
$$\begin{matrix}
H^i_{\mathrm{\acute{e}t}}(X^{\circ},\Q_{\ell})& \to & H^i(X')_{\ell}& \to & H^i(X)_{\ell}\\
\downarrow && \downarrow&& \downarrow\\
H^{i+2c}_{X^{\circ}}(Y^{\circ}, \Q_{\ell})(c)& \to & H^{i+2c}_{X'}(Y)_{\ell}(c)& \to & H^{i+2c}(Y)_{\ell}(c).
\end{matrix}$$
The left square is evidently commutative. 
  The composition 
$H^i_{\mathrm{\acute{e}t}}(X^{\circ}, \Q_{\ell})\to 
H^{i+2c}_{\mathrm{\acute{e}t}}(Y^{\circ}, \Q_{\ell})(c)\to H^{i+2c}(Y)_{\ell}(c)$ coincides with the composition 
$H^i_{\mathrm{\acute{e}t}}(X^{\circ}, \Q_{\ell})\to H^{i+2c}_{X^{\circ}}(Y^{\circ}, \Q_{\ell})(c)\to H^{i+2c}_{X'}(Y)_{\ell}(c)\to H^{i+2c}(Y)_{\ell}(c)$. Hence it is sufficient to prove the commutativity of the right square.
Let $p:=\dim(X)$, so $\dim(Y)=p+c$. 
  Let $j=2p-i$. It is sufficient to prove that for $a \in H^i(X')_{\ell}$ and $b\in H^j(Y)_{\ell}(p)$, we have
$(a \cup \xi \cup b)_Y=(a\cup b|_X)_X$ in $\Q_{\ell}$. Using $z=a\cup b|_{X'}\in H^{2p}(X')_{\ell}(p)$, we see that it is sufficient to prove that for $z\in H^{2p}(X')_{\ell}(p)$, the image of $z$ under $H^{2p}(X')_{\ell}(p)\to H^{2p+2c}_{X'}(Y)_{\ell}(p+c) \to H^{2p+2c}(Y)_{\ell}(p+c)\to \Q_{\ell}$ (the first arrow is the product with $\xi$) and the image of $z$ under $H^{2p}(X')_{\ell}(p)\to H^{2p}(X)_{\ell}(p)\to \Q_{\ell}$ coincide. 
 $H^{2p}(X')_{\ell}(p)$ is generated by the Chern classes of the $\cO_X$-modules $[\kappa(u)]$, where $u$ ranges over all non-singular closed points of $X$ and $\kappa(u)$ is the residue field at $u$. For $z= [\kappa(u)]$, the image of $z$ in $H^{2p+2c}(Y)_{\ell}(p+c)$ is the Chern class of the $\cO_Y$-module $\kappa(u)$. Hence the image of this $z$ in $\Q_{\ell}$ via $H^{2p+2c}(Y)_{\ell}(p+c)$ is $1$. On the other hand, the image of this $z$ in $\Q_{\ell}$ via $H^{2p}(X)_{\ell}(p)$ is $1$. 
  Thus both images coincide. 
\end{pf}

\begin{sbcor} Let $X$ be a projective vertical log smooth fs log scheme over a geometric standard log point $s$. Let $X'$ be 
a log blow-up of $X$ such that $M_{X'}/\cO_{X'}^\times$ is free ($\ref{moreg}$).  
Then 
the image of the Chern class map $\gr^i K_{\lim}(X) \to H^{2i}(X)_{\ell}(i)$ coincides with the image of the Chern class map $\gr^i K(X') \to H^{2i}(X)_{\ell}(i)$. 

\end{sbcor}

\begin{pf} Let $Y$ be any log blow-up of $X$ and let $a\in \gr^i K(Y)_\Q$. Take a log blow-up $Y'$ of $Y$ such that $M_{Y'}/\cO_{Y'}^\times$ is free and such that $Y'$ is also a log blow-up of $X'$. Let $a'$ be the image of $a$ in $\gr^i K(Y')$ by pullback, and let $b$ be the image of $a'$ in $\gr^i K(X')_\Q$ by pushforward. Then by Proposition \ref{RR}, the image of $a$ in $H^{2i}(X)_{\ell}(i)$ coincides with the image of $b$.\end{pf}

\begin{sbpara}\label{Tr} The above Proposition \ref{RR} contains the following trace formula in \cite{KatoSaito}. Let $X$ be a projective vertical log smooth fs log scheme over a geometric standard log point $s$. Assume that $X$ is purely of dimension $d$. Let $(X\times X)'$ be a log blow-up of $X\times X$, let $\alpha\in \gr^d K((X \times X)')_\Q$,  and let $f_{\alpha}$ be the image of $\alpha$ under the composition  $\gr^d K_{\lim}(X \times X) _\Q\to H^{2d}(X\times X)_{\ell}(d)\cong \bigoplus_i \; \Hom(H^i(X)_{\ell}, H^i(X)_{\ell})$, where the last isomorphism is by Poincar\'e duality (Corollary \ref{Q_lPoin}) and K\"unneth formula (Corollary \ref{Q_lKunn}). We consider the trace $\Tr(f_{\alpha})$. 
Let $X'$ be the log blow-up $X\times_{X\times X} (X\times X)'$ of the diagonal, and let the intersection of $\alpha$ with the diagonal $\alpha \cdot \Delta_X\in \Q$ be the image of $\alpha$ under the composition $\gr^d  K((X \times X)')_\Q\to \gr^d  K(X')_\Q \to K(s)_\Q=\Q$, where the first arrow is the pullback by $X' \to (X\times X)'$ and the second arrow is the pushforward. Then we have the trace formula
$$\Tr(f_{\alpha})= \alpha\cdot \Delta_X \in \Q.$$

This follows from Proposition \ref{RR} as follows. Consider the diagram 
$$\begin{matrix}  
&&\gr^d K_{\lim}(X \times X)_\Q & \to & \gr^d K_{\lim}(X)_\Q &\to &  \gr^0K(s)_\Q=\Q\\
&&\downarrow && \downarrow && \downarrow \\
\bigoplus_i \;\Hom(H^i(X)_{\ell}, H^i(X)_{\ell}) & \cong & H^{2d}(X \times X)_{\ell}(d)  & \to & H^{2d}(X)_{\ell}(d) & \to & H^0(s)_{\ell}=\Q_{\ell},
\end{matrix}$$
where the  first arrow in the lower row is the pullback by the diagonal. The left square is clearly commutative and the right square is commutative by Proposition \ref{RR}. The image of $f_{\alpha}\in \bigoplus_i \Hom(H^i(X)_{\ell}, H^i(X)_{\ell})$ in $\Q_{\ell}$ under the composition of the lower row is $\Tr(f_{\alpha})$. 
 This gives a proof of the trace formula. 

\end{sbpara}

\section{Log motives}
\label{s:logmot}

In this Section \ref{s:logmot}, 
let $S$ be an fs log scheme and let $\ell$ be a prime number which is invertible on $S$. 
  We define and study the category of log (pure) motives. 

\subsection{The category of log motives}
\label{LM}
We define the category of log motives over $S$. 

\begin{sbpara} 
  For a projective vertical log smooth fs log scheme $X$ over $S$ and for $r\in \Z$, consider the symbol $h(X)(r)$. 

Let $$h(X)(r)_{\ell}:=\bigoplus_m\; H^m(X)_{\ell}(r),\quad \mathrm{where}\quad H^m(X)_{\ell}= R^mf_*\Q_{\ell}\ (\mathrm{cf.\ }\ref{Chern})$$ with $f:X\to S$ and with $R^mf_*$ for the log \'etale topology. This is a smooth $\Q_{\ell}$-sheaf on the  log \'etale site of $S$ (see \ref{smoothQ_l}). 
\end{sbpara}
 
\begin{sbpara}
Let $X$ and $Y$ be projective vertical log smooth fs log schemes over a geometric standard log point (\ref{slp}).  
Let $r,s\in \Z$.

An element $\alpha$ of $\gr^iK_{\lim}(X\times Y)_\Q$ with $i=d+s-r$, where $d=\dim(X)$ induces a homomorphism $h(X)(r)_{\ell}\to h(Y)(s)_{\ell}$ as follows.

Let $\beta$ be the image of $\alpha$ under the Chern class map
$$\gr^iK_{\lim}(X\times Y)_{\Q}\to H^{2i}(X \times Y)_{\ell}(i).$$
Then for $m,n\in \Z$ such that $m-2r=n-2s$, we have 
the composition
$$H^m(X)_{\ell}(r) \to H^m(X \times Y)_{\ell}(r) \to H^{m+2i}(X \times Y)_{\ell}(r+i)$$
$$ \to H^{m+2i-2d}(Y)_{\ell}(r+i-d)= H^n(Y)_{\ell}(s).$$
Here the first arrow is the pullback, the second arrow is the cup product with $\beta$, the third arrow is the pushforward by the projection $X\times Y \to Y$. 
This gives a map $h(X)(r)_{\ell} \to h(Y)(s)_{\ell}$.
\end{sbpara}

\begin{sbpara}\label{11} 
Let $X$ and $Y$ be projective vertical log smooth fs log schemes over  $S$ and let $r,s\in \Z$. 

By definition, a morphism $f: h(X)(r) \to h(Y)(s)$  is a homomorphism $f: h(X)(r)_{\ell} \to h(Y)(s)_{\ell}$ of $\Q_{\ell}$-sheaves such that for any geometric standard log point $p$  over $S$, the pullback $h(X_p)(r)_{\ell} \to h(Y_p)(s)_{\ell}$ of $f$ is induced by an element of $\gr^{d+s-r} K_{\lim}(X_p\times_p Y_p)_{\Q}$ with $d=\dim(X_p)$ in the above way. 
\end{sbpara}

\begin{sbprop}
\label{iddia1}
  $(1)$ The identity morphism $h(X)(r)_{\ell} \to h(X)(r)_{\ell}$ is a morphism $h(X)(r)\to h(X)(r)$.

  $(2)$ More generally, for a morphism $Y\to X$ over $S$, the induced map $h(X)(r)_{\ell} \to h(Y)(r)_{\ell}$ is a morphism $h(X)(r)\to h(Y)(r)$.
\end{sbprop}

\begin{pf} We may and do assume that $S$ is a geometric standard log point $s$. Let $d$ be the dimension of $X$. 

We prove (1). Let $Z=X\times X$ (the fiber product over $S=s$) and consider the fan $\Sig:=\Sig_Z$ associated to $Z$ (\ref{case2}). 
  By Proposition \ref{cfan}, there is a finite subdivision $\Sig'\to \Sig$ such that $X':=X \times_{\Sig} \Sig'\to \Sig'$ and $Z':= Z \times_{\Sig} \Sig'\to \Sig'$ are strict. Hence the morphism $X'\to Z'$ is a strict closed immersion. Since a strict closed immersion between log smooth schemes is a regular immersion as is seen as in the classical case (cf.\ \cite{KatoLog} {\sc Proposition} (3.10)), this morphism $X'\to Z'$ is a regular immersion. Consider the $\cO_{Z'}$-module $\cO_{X'}$ and its class $[\cO_{X'}]\in \gr^d K(Z')_\Q$ with $d=\dim(X)$. By Poincar\'e duality (Corollary \ref{Q_lPoin}) and by K\"unneth formula 
(Corollary \ref{Q_lKunn}), this class induces the identity map $h(X)_{\ell}(r)\to h(X)_{\ell}(r)$. 

(2) follows from (1). The homomorphism
$h(X)_{\ell}(r)\to h(Y)_{\ell}(r)$ associated to $f$ is induced by an element of $\gr^dK_{\lim}(X\times Y)_\Q$ with $d=\dim(X)$ which is obtained from the above element of $\gr^dK_{\lim}(X\times X)_\Q$ giving the identity morphism, by pulling back by $1\times f$. 
\end{pf}

\begin{sbpara}\label{iddia2}  
The above Proposition \ref{iddia1} explains the reason why we must use $K_{\lim}$ (not just $K$) in the definition of morphism of the category of log motives. For a projective vertical log smooth fs log scheme $X$ over a geometric standard log point $s$, the diagonal  $X\to X\times X$ is usually not a regular immersion and cannot define an element of $K(X \times X)$. We need a log modification $Z\to X\times X$ to have an element of $K(Z)$ corresponding to the diagonal, which gives the identity morphism $h(X) \to h(X)$. 
\end{sbpara}

\begin{sbprop} For morphisms $f: h(X_1)(r_1) \to h(X_2)(r_2)$ and $g: h(X_2)(r_2)\to h(X_3)(r_3)$, the composition 
$g\circ f: h(X_1)(r_1)\to h(X_3)(r_3)$ is a morphism. 

\end{sbprop}

\begin{pf} 
We may assume that $S$ is a geometric standard log point. 
 If $f$ is induced by $\alpha\in \gr K_{\lim}(X_1\times X_2)_\Q$ and $g$ is induced by $\alpha'\in \gr K_{\lim}(X_2\times X_3)_\Q$, $g\circ f$ is induced by the following  element $\alpha''$ of $\gr K_{\lim}(X_1\times X_3)_\Q$. 
 Let 
$u\in \gr K_{\lim}(X_1\times X_2\times X_3)_\Q$ be the product of the pullbacks of $\alpha$ and $\alpha'$. 
Let $(X_1\times X_3)'$ be a log blow-up of $X_1\times X_3$ having free  $M/\cO^\times$ (\ref{moreg}), and let $(X_1\times X_2\times X_3)'$ be a log blow-up of $X_1\times X_2\times X_3$ having free 
$M/\cO^\times$ such that $u$ comes from an element $v$ of $\gr K((X_1\times X_2\times X_3)')_\Q$
and such that we have a morphism $(X_1\times X_2\times X_3)'\to (X_1\times X_3)'$ which is compatible with the projection $X_1\times X_2\times X_3\to X_1\times X_3$.
 Let $\alpha''$ be the
 pushforward of $v$ by the morphism $(X_1 \times X_2 \times X_3)'\to (X_1 \times X_3)'$. 
Then $g\circ f$ is induced by $\alpha''$ by Proposition \ref{RR}.
\end{pf}

\begin{sbpara}\label{deflmot}  Imitating the definition of motives by Grothendieck, we define {\it the category $\mathrm{LM}(S)$ of  log motives over $S$} as the category of the symbols $(h(X)(r), e)$, where $e$ is an idempotent in the endomorphism ring of $h(X)(r)$. The set of morphisms is defined as $$\Hom((h(X_1)(r_1), e_1), (h(X_2)(r_2), e_2)):= e_2\circ \Hom(h(X_1), h(X_2))\circ e_1\subset \Hom(h(X_1), h(X_2)).$$ The identity morphism of $(h(X), e)$ is $e$. 

The $\ell$-adic realization $M_{\ell}$ of the log motive $M=(h(X), e)$ is defined to be $e h(X)_{\ell}$. 
\end{sbpara}

\begin{sbpara} 
In the case where $S$ is $\Spec(k)$ for a field $k$ with trivial log structure, our category $\mathrm{LM}(S)$ coincides with the category of motives over $k$ modulo homological equivalence defined by Grothendieck. 

This is because $\CH^r(X \times Y)_{\Q} = \gr^r K(X\times Y)_\Q=\gr^rK_{\lim}(X\times Y)_\Q$.
\end{sbpara}

\begin{sbpara}
For a morphism $S'\to S$ of fs log schemes, we have the evident pullback
functor $\mathrm{LM}(S)\to \mathrm{LM}(S')$. 
\end{sbpara}

\begin{sbpara} 
  For a Galois extensino $p' \to p$ of standard log points, we have
$$\Hom_{\mathrm{LM}(p)}(h(X)(r), h(Y))(s)) \overset{\cong}\to \Hom_{\mathrm{LM}(p')}(h(X')(r), h(Y')(s))^G,$$
where $X'$ and $Y'$ are the base-changed objects from $X$ and $Y$, $r, s \in \bZ$, and $(-)^G$ denotes the $G$-invariant part for $G=\Gal(p'/p)$. 
\end{sbpara}

\subsection{Basic things}
\begin{sbpara}\label{dsum1} Direct sums and direct products exist in $\mathrm{LM}(S)$, and they coincide. 

In fact, we have $h(X) \oplus h(Y):= h(X \coprod Y)$, and 
if $r\leq s$, $h(X)(r)\oplus h(Y)(s) = (h((X\times {\bf P}^n)\coprod Y)(s), e)$ for $n\geq s-r$ and for some $e$. 
\end{sbpara}

\begin{sbconj}\label{conjoplus}
For a projective vertical log smooth fs log scheme $X$ of
relative dimension $d$ over $S$,
$h(X)$ has a decomposition
$$ h(X) = h^0(X) \oplus h^1(X) \oplus \dots \oplus h^{2d}(X) $$
in the category $\mathrm{LM}(S)$ of log motives such that $h^i(X)_{\ell}=H^i(X)_{\ell}$. 
\end{sbconj}

Note that such a decomposition is unique if it exists. 

\begin{sbpara} 
  We have the following: 
$h({\bf P}^n)= \bigoplus_{i=0}^n\;  h^{2i}({\bf P}^n)$. $h^{2i}({\bf P}^n)\cong \Q(-i)$ canonically for $0\leq i\leq n$. Here $\Q=h(S)$. 
\end{sbpara}

\begin{sbpara} We define the category $\mathrm{LM}^{\spl}(S)$ as follows. For a projective vertical log smooth fs log scheme $X$ over $S$ and for $m, r\in \Z$, consider the symbol $h^m(X)(r)$.

For projective vertical log smooth fs log schemes $X$ and $Y$ over $S$ and for $m,n, r,s\in \Z$, a morphism $h: h^m(X)(r)\to h^n(Y)(s)$ means a homomorphism $H^m(X)_{\ell}(r) \to H^n(Y)_{\ell}(s)$ of smooth $\Q_{\ell}$-sheaves on $S$ satisfying the following condition. If $m-2r\neq n-2s$, then $h=0$. If $m-2r=n-2s$, then for any geometric standard log point $p$ over $S$, the pullback of $h$ to $p$ comes from an element of $\gr^{d+s-r}K_{\lim}(X_p\times_p Y_p)$, where $d=\dim(X_0)$. %

  An object of $\mathrm{LM}^{\spl}(S)$ is $(h^m(X)(r), e)$, where $X$ is a projective vertical log smooth fs log scheme over $S$, $m, r\in \Z$, and $e$ is an idempotent of the ring of endomorphism of $h^m(X)(r)$. 
  Morphisms are defined like the case of $\mathrm{LM}(S)$.
\end{sbpara}

\begin{sbpara}\label{dsum2}  
  Similarly to the case of $\mathrm{LM}(S)$ (\ref{dsum1}), direct sums exist in $\mathrm{LM}^{\spl}(S)$.
  We have a functor $$\mathrm{LM}(S) \to \mathrm{LM}^{\spl}(S)\;;\; h(X)(r) \mapsto \bigoplus_m \; h^m(X)(r).$$ 

  Conjecture \ref{conjoplus} is that this functor is an equivalence of categories. 
\end{sbpara}

\begin{sbpara} Tensor products are defined in $\mathrm{LM}(S)$ as follows.

$(h(X)(r), e) \otimes (h(X')(s),e'):= (h(X \times X')(r+s), e\otimes e')$.

For a log motive $M$ over $S$, the Tate twist $M(-r)$ ($r\geq 0$) is identified with 
$M\otimes h^{2r}({\bf P}^n)$ with $n\geq r$. 

\end{sbpara}

\begin{sbpara}
  Compared with $\mathrm{LM}(S)$, a disadvantage of the category $\mathrm{LM}^{\spl}(S)$ is that the tensor products cannot be defined. 
\end{sbpara}

\begin{sbpara}
  Duals are defined in $\mathrm{LM}(S)$ as follows. 

$(h(X)(r),e)^*=(h(X)(d-r),e(d-2r))$, where $d$ is the relative dimension of $X$ over $S$.
  
  Note that any morphism $h(X)(r) \to h(Y)(s)$ induces a homomorphism 
$(h(Y)(s)^*)_{\ell} \to (h(X)(r)^*)_{\ell}$ of $\bQ_{\ell}$-sheaves by Poincar\'e duality (Corollary \ref{Q_lPoin}). 
  We can easily check that this homomorphism 
gives a morphism $h(Y)(s)^*=h(Y)(d'-s) \to h(X)(d-r)=h(X)(r)^*$ of motives, 
where $d'$ is the relative dimension of $Y$ over $S$ 
by using the same elements of $\gr^iK_{\lim}(X_p \times_p Y_p)_{\bQ}$, 
where $p$ is a geometric standard log point over $S$ and $i=d+s-r=d'+(d-r)-(d'-s).$
\end{sbpara}

\begin{sbpara}\label{MW1}  
  Let $X$ be a projective vertical log smooth fs log scheme over $S$.
  We conjecture that, for any morphism 
$s \to S$ from a standard log point associated to some finite field and 
for each $m\in \bZ$, the filtration (the monodromy filtration) on the stalk over $s$ of $H^m(X)_{\ell}$ 
determined by the monodromy operator coincides with the Frobenius weight filtration. 
  We call this the {\it monodromy-weight conjecture} for $X$.
\end{sbpara}

\begin{sbprop}\label{MW2} 
  Let $X$ and $Y$ be projective vertical log smooth fs log schemes over $S$. 
  Assuming the monodromy-weight conjecture for $X$ and $Y$, %
we have the following$:$

\noindent 
If $m-2r>n-2s$ and if $S$ is of finite type over $\Z$, there is no non-zero homomorphism $H^m(X)_{\ell}(r)\to H^n(Y)_{\ell}(s)$. 
\end{sbprop}

\begin{pf} 
  This is reduced to the case where $S$ is a standard log point associated to a finite field $k$. 
  Let $w=m-2r, w'=n-2s$. 
  Monodromy-weight conjecture tells that as a finite-dimensional $\Q_{\ell}$-vector space  with actions of $\Gal(\bar k/k) $ and the monodromy operator $\cN$, 
the stalk of $H^m(X)_{\ell}(r)$ (resp.\ $H^n(Y)_{\ell}(s)$) 
is isomorphic to a direct sum of subobjects $Q$ (resp.\  $R$) being isomorphic to $\mathrm{Sym}^i H^1(E)_{\ell}\otimes V$, where the action of $\Gal(\bar k/k)$ on $V$ is of weight $w-i$ (resp.\  $w'-i$) and the action of $\cN$ on $V$ is trivial. 
  Hence, it is enough to show that there is no non-zero $\Q_{\ell}$-linear map $Q \to R$ which is compatible with the actions of $\Gal(\bar k/k)$ and $\cN$. 
  Let $Q \to R$ be such a map. 
  For any non-zero element $x\in R$ of weight $u\geq w'$, we have $\cN^{u-w'}(x)\neq 0$. But as a $\Q_{\ell}$-vector space with an action of $\cN$, $Q$ is generated by an element $y$ of weight $u\geq w$ such that $\cN^{u-w+1}(y)=0$. The image $x$ of this $y$ in $R$ is of weight $u \geq w'$ and $\cN^{u-w'}(x)=0$ because $u-w'\geq u-w+1$. Hence $x=0$. 
  Therefore the map $Q\to R$ is the zero map. 
\end{pf}

\begin{sbrem}
\label{nontrivhom}
  On the other hand, a non-trivial homomorphism $H^m(X)_{\ell}(r) \to H^n(Y)_{\ell}(s)$ can exist even if $m-2r<n-2s$
and even if $S$ is of finite type over $\bZ$. 
  In fact, let $S$ be a standard log point, $X=S$, and $Y$ the log Tate curve. 
  Then we have an exact sequence $0 \to \bQ_{\ell} \to H^1(Y)_{\ell} \to \bQ_{\ell}(-1) \to 0$. 
  Hence a non-trivial homomorphism $H^0(X)_{\ell}\to H^1(Y)_{\ell}$ exists. 
\end{sbrem}

\begin{sbpara}\label{logpic0}
For an $X$ semistable over a standard log point, $H^1(X_{\mathrm{Zar}}, M^{\gp}_X/\cO^\times_X)=0$ because $M^{\gp}_X/\cO^\times_X$ is a flasque sheaf, which implies that $\mathrm{Pic}(X)=H^1(X_{\mathrm{Zar}}, \cO_X^\times) \to H^1(X_{\mathrm{Zar}}, M_X^{\gp})$ is surjective. Hence by \ref{TS}, we have:

Let $X, Y$ be  projective vertical log smooth fs log schemes over an fs log scheme $S$. Then an element of $H^1((X\times Y)_{\mathrm{Zar}}, M^{\gp}_{X \times Y})$ gives a homomorphism $h(X)(r) \to h(Y)(r+1-d)$, where $r\in \Z$ and $d$ is the relative dimension of $X$ over $S$.  
\end{sbpara}

\subsection{Numerical equivalence}
\begin{sbprop}
\label{TrQ}
 For any log motive $M$ over $S$ and for any 
morphism $f: M\to M$, $\Tr(f)\in\Q_{\ell}$ belongs to $ \Q$. 
(Precisely speaking, $\Tr(f)$ is a locally constant function  $S\to \Q$. It is constant if $S$ is connected.)
\end{sbprop}

\begin{pf}
We are reduced to the case where $S$ is a geometric standard  log point. Then the result follows from the trace formula \ref{Tr}.
\end{pf}

\begin{sbpara}\label{nmeq1}  Definition of numerical equivalence. 

For objects $M$ and $M'$ of $\mathrm{LM}(S)$ and for a morphism $f: M\to M'$, we say $f$ is {\it numerically equivalent to $0$} if for any morphism $g: M'\to M$, we have $\Tr(gf)=0$, that is, $\Tr(fg)=0$. 

Morphisms $f, g:M\to M'$ are said to be {\it numerically equivalent} if $f-g$ is numerically equivalent to $0$. 
\end{sbpara}

\begin{sblem}\label{nmeq2} 
  Let $\sim$ be the numerical equivalence. Let $f,g:M\to N$ be morphisms in $\mathrm{LM}(S)$. Then 

$(1)$ $fh\sim gh$ for any morphism $h: L\to M$ from a log motive $L$ over $S$.

$(2)$ $hf\sim hg$ for any morphism $h:N\to L$ to a log motive $L$ over $S$. 
\end{sblem}

This is easy to prove.

\begin{sbpara}
By Lemma \ref{nmeq2}, we have the category $\mathrm{LM}_{\mathrm{num}}(S)$ of log motives over $S$ modulo numerical equivalence. 
\end{sbpara}

\begin{sbconj}\label{n=h} 
  In $\mathrm{LM}(S)$, $f\sim g$ implies $f=g$. 
  That is, $\mathrm{LM}(S)= \mathrm{LM}_{\mathrm{num}}(S)$. 
\end{sbconj}

\begin{sbpara}
  When $S$ is a geometric standard log point, 
the category $\mathrm{LM}_{\mathrm{num}}(S)$ is independent of the choice of $\ell$.
  This is a consequence of Proposition \ref{TrQ} since in this case, the group $\Hom(h(X)(r),h(Y)(s))$ is identified with a quotient of $\gr^{d+s-r}K_{\lim}(X \times_SY)_{\bQ}$ in the notation in \ref{11}. 
\end{sbpara}

\subsection{Semisimplicity}
\label{semisimple}
\begin{sbthm}
\label{thmss} 
$(1)$ The category $\mathrm{LM}_{\mathrm{num}}(S)$  is a semisimple abelian category. 

\medskip

$(2)$ The category $\mathrm{LM}(S)$   is a semisimple abelian category if and only if the numerical equivalence for morphisms of this category is trivial. 

\end{sbthm}
To prove this, 
we imitate the method of U.\ Jannsen in \cite{Jannsen}.

\begin{sbpara} The following fact is known. %

A pseudo-abelian category $\cC$ is a semisimple abelian category if the following (i) and (ii) are satisfied for any objects $X$ and $Y$.

(i) $\Hom\,(X,Y)$ is a $\Q$-vector space, the composition of morphisms is bilinear, and any idempotent of $\End\,(X)$ has a kernel.

(ii) $\End\,(X)$ is a finite-dimensional semisimple $\Q$-algebra. 
\end{sbpara}

\begin{sblem}\label{sslem1} Let $F$ be a field, $A$, $B$ finite-dimensional $F$-vector spaces, $(\;,\;): A\times B \to F$ an $F$-bilinear map, $F_0$ a subfield of $F$, $A_0$ an $F_0$-subspace of $A$, and $B_0$ an $F_0$-subspace of $B$. Assume that $A$ is generated by $A_0$ over $F$,  $B$ is generated by $B_0$ over $F$, and $(a,b)\in F_0$ for any $a\in A_0$ and $b\in B_0$. Let 
 $K=\{a\in A\; |\; (a,b)=0\;\text{for any} \; b\in B\}$, and $K_0= \{a\in A_0\;|\; (a, b)=0 \;\text{for any}\; b\in B_0\}$. Then$:$

\medskip

 $F\otimes_{F_0} A_0/K_0 \overset{\cong}\to A/K$. 

In particular, $A_0/K_0$ is finite-dimensional over $F_0$. 
\end{sblem}

\begin{pf} Take an $F_0$-subspace $A_0'$ of $A_0$ such that $F\otimes_{F_0} A_0'\overset{\cong}\to A$ and an $F_0$-subspace $B_0'$ of $B_0$ such that $F\otimes_{F_0} B_0' \overset{\cong}\to B$. Then $A_0'$ and $B'_0$ are finite-dimensional over $F_0$. Let $K'_0= \{a\in A_0'\:|\; (a,b)=0\; \text{for any} \; b\in B\}=  \{a\in A_0'\:|\; (a,b)=0\; \text{for any} \; b\in B'_0\}$. Let 
$L'_0=\{b \in B_0'\; |\; (a,b)=0\;\text{for any}\; a\in A\}=\{b \in B_0'\; |\; (a,b)=0\;\text{for any}\; a\in A_0'\}$. The composition
$$A_0'/K_0' \to A_0/K_0 \to \Hom(B_0'/L_0', F_0)$$
is an isomorphism and the two arrows here are injective. Hence 
we have

\medskip

$(*)$ $A'_0/K_0'\to A_0/K_0$ is an isomorphism. 

\medskip

On the other hand, the paring $A \times B \to F$ is identified with $F\otimes_{F_0}$ of the pairing $A_0'\times B'_0\to F_0$. Hence we have

\medskip

$(**)$ $F \otimes_{F_0} A'_0/K'_0 \to A/K$ is an isomorphism.

\medskip

By $(*)$ and $(**)$, we have that $F\otimes_{F_0} A_0/K_0 \to A/K$ is an isomorphism. 
\end{pf}

\begin{sblem}\label{sslem2} Let $F$ be a field of characteristic $0$, 
 $V$ a finite-dimensional $F$-vector space, and $A$ an $F$-subalgebra of $\End_F(V)$. Let $J$ be the Jacobson radical of $A$, that is, $J$ is the largest nilpotent two-sided ideal of $A$. Let $I= \{a\in A\;|\; \Tr(ab)=0\;\text{for any}\; b\in A\}$. Here $\Tr$ is the trace of an $F$-linear map $V\to V$. Then $I=J$.

\end{sblem}

\begin{pf} Let $a\in J$. Then for any $b\in A$, $ab$ is nilpotent and hence $\Tr(ab)=0$. Hence $a\in I$. Next we prove $I\subset J$. We may assume that $F$ is algebraically closed. It is sufficient to prove that all elements of $I$ are nilpotent. Let $a\in I$. Let $(\alpha_i)_{1\leq i \leq n}$ ($n=\dim_F(V)$) be the eigenvalues of $a$ counted with multiplicity. We have $0=\Tr(a^n)=\sum_{i=1}^n \alpha_i^n$ for any $n\geq 1$. This proves that $\alpha_i=0$ for all $i$. Hence $a$ is nilpotent. 
\end{pf}

\begin{sblem}\label{sslem3} Let $F$ be a field of characteristic $0$, 
 $V$ a finite-dimensional $F$-vector space, $A$ an $F$-subalgebra of $\End_F(V)$, $F_0$ a subfield of $F$, and $A_0$ an $F_0$-subalgebra of $A$. Assume that $A_0$ generates the $F$-vector space $A$ and assume that $\Tr(a)\in F_0$ for any $a\in A_0$. Let $I_0=\{a\in A_0\;|\; \Tr(ab)=0\;\text{for any}\; b\in A_0\}$. Then $I_0$ is a two-sided ideal of $A_0$, $A_0/I_0$ is a finite-dimensional semisimple $F_0$-algebra, and all elements of $I_0$ are nilpotent. 

\end{sblem}

\begin{pf} The fact that $I_0$ is a two-sided ideal of $A_0$ is shown easily. Let $I=\{a\in A\;|\; \Tr(ab)=0\;\text{for any}\; b\in A\}$. Then $I$ is nilpotent  and $A/I$ is a semisimple algebra by Lemma \ref{sslem2}. Hence all elements of $I_0$ are nilpotent. 
By  Lemma \ref{sslem1}, $A_0/I_0$ is finite-dimensional and $F\otimes_{F_0} A_0/I_0$ is isomorphic to $A/I$. Hence $A_0/I_0$ is semisimple. 
\end{pf}

\begin{sbpara} We prove Theorem \ref{thmss} (1). Let $M$ be a log motive over $S$. In Lemma \ref{sslem3}, take $F=\Q_{\ell}$, $F_0=\Q$, and let  $A$ be the $\Q_{\ell}$-subalgebra of $\mathrm{End}_{\Q_{\ell}}(M_{\ell})$ generated by $A_0:= \mathrm{End}_{\mathrm{LM}(S)}(M)$. Then the endomorphism ring of $M$ in the category of log motives over $S$ modulo numerical equivalence is $A/I_0$, where $I_0$ is as in Lemma \ref{sslem3}. By Lemma \ref{sslem3}, $A/I_0$ is a finite-dimensional semisimple $\Q$-algebra. This proves (1) of Theorem \ref{thmss}.

We prove Theorem \ref{thmss} (2). The if part follows from (1). We prove the only if part. Let $F=\Q_{\ell}$, $F_0=\Q$, and $ A, A_0, I_0$ be as in the proof of (1). By Lemma \ref{sslem3}, all elements of  $I_0$ are nilpotent. Assume that $A_0$ is semisimple. Since $I_0$ is a two-sided ideal of $A_0$ and all elements of $I_0$ are nilpotent, we have $I_0=0$. That is, the numerical equivalence is trivial.

\end{sbpara}

\section{Log mixed motives}
\label{s:logmm}

  We define the category of log mixed motives. 

\subsection{The category $\cC_S$}

\begin{sbpara}
\label{categoryC_S}
  Let $\ell$ be a prime number. 
  Let $S$ be an fs log scheme over $\Z[1/\ell]$ of finite type. 

  Let $\cC_S$ be the following category.

  Objects:  $(\cF, W, (X_w)_{w\in \Z},  (V_{w,1})_{w\in \Z}, (V_{w,2})_{w\in \Z}, (\iota_w)_{w\in \Z})$. 

  Here $\cF$ is a smooth $\Q_{\ell}$-sheaf on the log \'etale site of $S$. The $W$ is an increasing filtration on $\cF$ by smooth $\Q_{\ell}$-subsheaves. The $X_w$ is a projective vertical log smooth fs log scheme over $S$. For each $w\in \Z$, $V_{w,1}$ and $V_{w,2}$ are smooth $\Q_{\ell}$-subsheaves of $\bigoplus_{r\in \Z} 
H^{w+2r}(X_w)_{\ell}(r)$ such that $V_{w,1}\subset V_{w,2}$. 
  The $\iota_w$ is an isomorphism $\gr^W_w\cF\cong V_{w,2}/V_{w,1}$. 

$W$ is called the {\it weight filtration}. 

A morphism $$(\cF, W, (X_w)_{w\in \Z}, (V_{w,1})_{w\in \Z}, (V_{w,2})_{w\in \Z}, (\iota_w)_{w\in \Z})\to (\cF', W', (X'_w)_{w\in \Z}, (V'_{w,1})_{w\in \Z}, (V'_{w,2})_{w\in \Z}, (\iota'_w)_{w\in \Z})$$ in $\cC_S$  is a homomorphism of $\Q_{\ell}$-sheaves $\cF\to \cF'$ which respects the weight filtrations such that for each $w\in \Z$, the pullback of $\gr^W_w\cF\to \gr^{W'}_w\cF'$ to any geometric standard log point $s$ over $S$ is induced from the sum of morphisms 
$h(X_w\times_S s)(r)\to h(X'_w \times_S s)(r')$ for various $r, r' \in \bZ$ which sends $V_{w,i}$ to $V'_{w,i}$ over $s$ for $i=1,2$. 
\end{sbpara}

\begin{sbpara}
  The category $\cC_S$ has $\oplus$, the kernels and the cokernels. 

  Furthermore, $\otimes$, the dual, and Tate twists are defined in $\cC_S$. 

  These are explained in \ref{cC1}--\ref{cC5}. 
\end{sbpara}

\begin{sbpara}\label{cC1}
We have 
$$(\cF, W, (X_w)_{w\in \Z}, (V_{w,1})_{w\in \Z}, (V_{w,2})_{w\in \Z}, (\iota_w)_{w\in \Z})\oplus(\cF', W', (X'_w)_{w\in \Z}, (V'_{w,1})_{w\in \Z}, (V'_{w,2})_{w\in \Z}, (\iota'_w)_{w\in \Z})$$ 
$$=(\cF\oplus \cF', W\oplus W', (X_w\coprod X'_w)_{w\in \Z}, (V_{w,1}\oplus V'_{w,1})_{w\in \Z}, (V_{w,2}\oplus V'_{w,2})_{w\in \Z}, (\iota_w\oplus \iota'_w)_{w\in \Z}).$$
\end{sbpara}

\begin{sbpara}\label{cC2}
The kernel of 
$$(\cF, W, (X_w)_{w\in \Z}, (V_{w,1})_{w\in \Z}, (V_{w,2})_{w\in \Z}, (\iota_w)_{w\in \Z})\to (\cF', W', (X'_w)_{w\in \Z}, (V'_{w,1})_{w\in \Z}, (V'_{w,2})_{w\in \Z}, (\iota'_w)_{w\in \Z})$$
is $(\cF'', W'', (X''_w)_{w\in \Z}, (V''_{w,1})_{w\in \Z}, (V''_{w,2})_{w\in \Z}, (\iota''_w)_{w\in \Z})$,
where $\cF''$ is the kernel of $\cF\to \cF'$, $W''$ is induced from $W$, $X''_w=X_w$, $V''_{w,2}$ is the kernel of $V_{w,2}\to V'_{w,2}/V'_{w,1}$, $V''_{w,1}= V''_{w,2}\cap V_{w,1}$, and $\iota_w''$ is induced from $\iota_w$. 
\end{sbpara}

\begin{sbpara}\label{cC3}
The cokernel of the above morphism is $(\cF'', W'', (X''_w)_{w\in \Z}, (V''_{w,1})_{w\in \Z}, (V''_{w,2})_{w\in \Z}, (\iota''_w)_{w\in \Z})$,
where $\cF''$ is the cokernel of $\cF\to \cF'$, $W''$ is induced from $W'$, $X''_w=X'_w$, $V''_{w,2}= V'_{w,2}+\mathrm{Image}(V_{w,2})$, $V''_{w,1}=V'_{w,1}+\mathrm{Image}(V_{w,2})$, and 
$\iota_w''$ is induced by $\iota'_w$. 
\end{sbpara}

\begin{sbpara}\label{cC4}
$$(\cF, W, (X_w)_{w\in \Z}, (V_{w,1})_{w\in \Z}, (V_{w,2})_{w\in \Z}, (\iota_w)_{w\in \Z})\otimes (\cF', W', (X'_w)_{w\in \Z}, (V'_{w,1})_{w\in \Z}, (V'_{w,2})_{w\in \Z}, (\iota'_w)_{w\in \Z})$$
is defined as $(\cF'', W'', (X''_w)_{w\in \Z}, (V''_{w,1})_{w\in \Z}, (V''_{w,2})_{w\in \Z}, (\iota''_w)_{w\in \Z})$, where $\cF''=\cF\otimes \cF'$, $W''$ is the convolution of $W$ and $W'$, $X''_w=\coprod_{i+j=w} X_i\times X'_j$, $V''_{w,2}= \bigoplus_{i+j=w} V_{i,2}\otimes  V'_{j,2}$, $V''_{w,1}= \bigoplus_{i+j=w} (V_{i,1}\otimes V'_{j,2}+ V_{i,2}\otimes V'_{j,1})$, $\iota''_w=\bigoplus_{i+j=w} \iota_i \otimes \iota'_j$. 
\end{sbpara}

\begin{sbpara}\label{cC5}
  The definition of the dual and the Tate twists are the evident ones.
\end{sbpara}

\subsection{The category of log mixed motives}
\label{ss:LMM}
In his papers \cite{D1} and \cite{D2}, Deligne showed how we can obtain mixed Hodge structures of geometric origin basing on the theory of pure Hodge structures.  
We imitate his method to formulate objects of $\cC_S$ of geometric origin. 

In this \ref{ss:LMM}, $S$ denotes an fs log scheme and $\ell$ denotes a prime number which is invertible on $S$. 

For an fs log scheme $X$ over $S$, $H^m(X)(r)_{\ell}$ denotes $R^mf_*\Q_{\ell}(r)$, where $f$ is the morphism $X\to S$. 

\begin{sbpara}\label{mix1}
Consider $(U, X, D)$, where $X$ is a projective vertical log smooth fs log scheme over $S$, $D=(D_\lam)_{\lam\in \Lam}$ is a finite family of  Cartier divisors on $X$, and $U$ is the open subscheme of $X$ defined as the complement of $\bigcup_{\lambda\in \Lam} D_\lam$ in $X$ satisfying the following condition:

For any subset $\Lam'$ of $\Lam$, $D_{\Lam'}:=\bigcap_{\lam\in \Lam'} D_\lam$ with the inverse image of the log structure of $X$ is log smooth over $S$, and of codimension $\sharp(\Lam')$ in $X$ at each point of it. 
\end{sbpara}

\begin{sbpara} Let the notation and the assumptions be as in \ref{mix1}. 

For $i\geq 0$, let $D^{(i)}$ be the disjoint union of $D_{\Lam'}$ for all $\Lam'\subset \Lam$ such that $\sharp(\Lam')=i$. 
In particular, $D^{(0)}=X$. 

For $i\geq 0$, we have a smooth $\Q_{\ell}$-sheaf $H^m(D^{(i)})_{\ell}$ on the  log \'etale site of $S$ (cf.\ \ref{smoothQ_l}). 
\end{sbpara}

\begin{sbpara}\label{mix2} 
  Let the notation and the assumptions be as in \ref{mix1}. Endow $U$ with the inverse image of the log structure of $X$. 

  Then $H^m(U)_{\ell}$ is a smooth $\Q_{\ell}$-sheaf on the log \'etale site of $S$ and 
we have a spectral sequence 
$$E_1^{i,j}=  H^{2i+j}(D^{(-i)})_{\ell}(i) \Rightarrow E_{\infty}^m= H^m(U)_{\ell}$$
in the category of smooth $\bQ_{\ell}$-sheaves. 
  In fact, first, by relative purity in log \'etale cohomology (\cite{HigashiyamaKamiya}), we have a spectral sequence 
with finite coefficients. 
  By Proposition \ref{prsm}, the $E_1$-terms of this spectral sequence determines a smooth $\bQ_{\ell}$-sheaves, which 
implies the above facts. 
\end{sbpara}

\begin{sbpara}\label{mix8} Consider a simplicial system $(U_{\bullet}, X_{\bullet}, D_{\bullet})$ of objects $(U, X, D)$ of \ref{mix1} (here we follow \cite{D2}). Let ${\bf H}^m(U_{\bullet})_{\ell}$ be the smooth $\Q_{\ell}$-sheaf on $S$ defined to be the $m$-th hypercohomology (relative to $S$) of the simplicial system. 
The  spectral sequence in \ref{mix2} is generalized to the spectral sequence 
$$E_1^{i,j}=  \bigoplus_{s\geq 0}\;  H^{j-2s}(D_{s+i}^{(s)})_{\ell}(-s) \Rightarrow E_{\infty}^m= {\bf H}^m(U_{\bullet})_{\ell}.$$
\end{sbpara}

\begin{sbpara}\label{mix3} Let the notation be as in \ref{mix8}. Let $m\in \Z$. We define an increasing filtration $W$ on ${\bf H}^m(U_{\bullet})_{\ell}$, which we call the {\it weight filtration}, as the filtration defined by the spectral sequence in \ref{mix8}. 
\end{sbpara}

\begin{sbpara}\label{mix4} If $S$ is of finite type over $\Z[1/\ell]$, let $\cC_S^{\mathrm{mot}}$ be the full subcategory of $\cC_S$ consisting of objects
which are obtained from the following standard objects in \ref{mix5} below by taking $\oplus$, kernels, cokernels, $\otimes$, the duals, and Tate twists. 

\end{sbpara}

\begin{sbpara}\label{mix5} 
  In the above, a standard object means:

  Consider $(U_{\bullet}, X_{\bullet}, D_{\bullet}, m)$, where $(U_{\bullet}, X_{\bullet}, D_{\bullet}) $ is as in \ref{mix8} and $m\in \Z$. 
  The associated standard object is as follows: 

  Let $\cF={\bf H}^m(U_{\bullet})_{\ell}$ on $S$. 

  Let $W$ be the filtration on ${\bf H}^m(U_{\bullet})_{\ell}$ defined by the spectral sequence in \ref{mix3}. 
  Then, for $w\in \Z$, $\gr^W_w{\bf H}^m(U_{\bullet})_{\ell}=V'_{w,2}/V'_{w,1}$ for some  $\Q_{\ell}$-subsheaves 
$V'_{w,1}$, $V'_{w,2}$ of $\bigoplus_{s\geq 0} H^{w-2s}(D_{s+m-w}^{(s)})_{\ell}(-s)$ such that 
$V'_{w,1} \subset V'_{w,2}$.

  Let $X_w=\bigsqcup_{s \geq 0} D^{(s)}_{s+m-w}$.

  Consider the natural projection 
$\bigoplus_{r\in \Z} H^{w+2r}(X_w)_{\ell}(r)\to
\bigoplus_{s\geq 0} H^{w-2s}(D_{s+m-w}^{(s)})_{\ell}(-s)$.  

  For $i=1, 2$, let $V_{w,i}$ be the pullbacks of $V'_{w,i}$ by this natural projection. 
  Then we have the isomorphism 
$\iota_w\colon \mathrm{gr}_w^W\cF \cong V_{w,2}/V_{w,1}$. 
\end{sbpara}

\begin{sbpara}\label{Cmot}
  If $S$ is affine and is the inverse limit of $S_{\lambda}$ which are of finite type over $\Z[1/\ell]$, we define $\cC_S^{\mathrm{mot}}$ as the inductive limit of the categories $\cC_{S_{\lambda}}^{\mathrm{mot}}$. 
  This does not depend on the choice of limits. 
\end{sbpara}

\begin{sbpara} 
  We define the category of log mixed motives $\mathrm{LMM}(S)$ over $S$ as the Zariski sheafification of the categories $\cC^{\mathrm{mot}}_S$ in \ref{Cmot}.
\end{sbpara}

\begin{sbpara} 
  For a morphism $S'\to S$ of fs log schemes, we have the pullback functor $\mathrm{LMM}(S) \to \mathrm{LMM}(S')$. 
\end{sbpara}

\begin{sbpara} We have a fully faithful functor  $$\mathrm{LM}^{\spl}(S) \to \mathrm{LMM}(S)$$
which sends $H^m(X)(r)$ to the object associated to $(U_{\bullet}, X_{\bullet}, D_{\bullet}, m)(r)$ with $X_{\bullet}$ determined by $X$, $U_{\bullet}= X_{\bullet}$, $D_{\bullet}$ empty.
\end{sbpara}

\begin{sbpara}\label{MMS} If the log structure of $S$ is trivial, we define the category $MM(S)$ of mixed motives to be the category of log mixed motives over $S$. 
\end{sbpara}

\subsection{Justifications of our definition}
\label{justification}

Here we explain the reason why we think our definition of log mixed motives is reasonable. 

\begin{sbpara} The reader may feel strange that 
in our definition of a morphism of log mixed motives (\ref{categoryC_S}), we do not put much conditions other than the condition that its $\gr^W$ is motivic, 
though it is usually impossible to take care of mixed objects by using only pure objects. 

We hope that the following Proposition \ref{justi0} (resp.\  Proposition \ref{justi})  justifies our definition of log mixed motive (resp.\  of morphism of log mixed motives) in \ref{ss:LMM} (resp.\ \ref{categoryC_S}). 

We hope that if $S$ is of finite type over $\Z$ and if we take the category of log mixed motives over $S$ and the category of smooth $\Q_{\ell}$-sheaves on the  log \'etale site of $S$ as $\cC_1$ and $\cC_2$, respectively, the conditions in \ref{C1C2} below are satisfied. 
  (Especially we hope that the finiteness assumption on $S$ assures that the condition (v) in \ref{C1C2} below is satisfied.)
\end{sbpara}

\begin{sbpara}\label{C1C2}
Let $\cC_1$ and $\cC_2$ be abelian categories.
 Assume that we have exact subfunctors $W_w: \cC_1\to \cC_1$ ($w\in \Z$) of the identify functor $\cC_1\to \cC_1$ such that $W_w\circ W_w=W_w$ and such that $W_{w'} \subset W_w$ if $w'\leq w$.  Assume that we have a functor $F: \cC_1\to \cC_2$. Assume that these satisfy the following six conditions. 

(i) For each object $M$ of $\cC_1$, $W_wM=M$ if $w\gg 0$ and $W_wM=0$ if $w\ll 0$.

(ii) The functor $F$ is exact.

(iii) Let $w\in \Z$ and let $M$ and $N$ be objects of $\cC_1$. Assume that $M$ and $N$ are  pure of weight $w$ (that is, $W_wM=M$, $W_{w-1}M=0$, $W_wN=N$, $W_{w-1}N=0$). Then the canonical map $\Hom_{\cC_1}(M, N) \to \Hom_{\cC_2}(F(M), F(N))$ 
is injective.

(iv) Let $w, w'\in \Z$ and assume $w>w'$. Let $M$ and $N$ be objects of $\cC_1$ and assume that $M$ is pure of weight $w$ and 
$N$ is pure of weight $w'$.  Then $\Hom(M, N)=0$ and $\Hom(F(M), F(N))=0$. 

(v) Let $w, w'\in \Z$ and assume $w\geq w'$. Let $M$ and $N$ be objects of $\cC_1$ and assume that $M$ is pure of weight $w$ and $N$ is pure of weight $w'$. Then the canonical map $\Ext^1_{\cC_1}(M, N) \to \Ext^1_{\cC_2}(F(M), F(N))$ is injective. 

(vi) Let $w\in \Z$. Then the full subcategory of $\cC_1$ consisting of all objects  which are pure of weight $w$ is semisimple.

\medskip

\noindent 
{\it Remark.}
  By Proposition \ref{MW2}, $\Hom(F(M), F(N))=0$ in the condition (iv) is reasonable. (This is clearly reasonable if the log structure of $S$ is trivial, but not trivial otherwise.) 
  Further, the condition (v) is related to Tate conjecture. 
  In fact, it means that an extension of motives splits if the $\ell$-adic realization splits; 
two extensions are isomorphic if their $\ell$-adic realizations are isomorphic. 
  These are Tate conjectures.
\end{sbpara}

\begin{sblem}\label{Claim0}
Let the notation and the assumptions be as in $\ref{C1C2}$ and let $M$ and $N$ be objects of $\cC_1$. 

$(1)$ The morphism $\Hom_{\cC_1}(M, N) \to \Hom_{\cC_2}(F(M), F(N))$ is injective. 

$(2)$ If there is a $w\in \Z$ such that $W_wM=0$ and $W_wN=N$, then $\Hom_{\cC_1}(M, N)=0$ and the map $\Ext^1_{\cC_1}(M, N) \to \Ext^1_{\cC_2}(F(M), F(N))$ is injective. 
\end{sblem}

This is easy to prove.

\begin{sbprop}\label{justi0}
Let the notation and the assumptions be as in $\ref{C1C2}$. Let $M$ be an object of $\cC_1$ and let $V$ be a subobject of $F(M)$ in $\cC_2$ such that for any $w\in \Z$, the subobject $\gr^W_wV:=(V\cap F(W_wM))/(V\cap F(W_{w-1}M))$ of  $F(\gr^W_wM)$ is $F(N_w)$ for some subobject $N_w$ of $\gr^W_wM$ in $\cC_1$. Then there is a unique subobject $N$ of $M$ in $\cC_1$ such that $V$ coincides with $F(N)$. 
\end{sbprop}

\begin{pf}  By downward induction on $w$, we may assume that $W_{w-1}M=0$  and that if we denote $V\cap F(W_wM)$ by $V'$, the subobject $V'':=V/V'$ of $F(M/W_wM)$ coincides with $F(N'')$ for some subobject $N''$ of $M/W_wM$. 
By the assumption, the subobject $V'$ of $F(W_wM)=F(\gr^W_wM)$ coincides with $F(N')$ for some subobject $N'$ of $W_wM$. 

Let the exact sequence $0 \to W_wM  \to U \to N'' \to 0$  be the pullback of the exact sequence 
$0\to W_wM \to M \to M/W_wM \to 0$  by $N'' \to M/W_wM$. Then $\mathrm{class}(F(U))\in \Ext^1_{\cC_2}(F(N''), F(W_wM))$ coincides with the image of $\mathrm{class}(V)\in \Ext^1_{\cC_2}(F(N''), F(N'))$ under the homomorphism $\Ext^1_{\cC_2}(F(N''), F(N')) \to \Ext^1_{\cC_2}(F(N''), F(W_wM))$ induced by the morphism $N' \to W_wM$. 
\medskip

{\bf Claim 1.} There are an object $N$ of $\cC_1$ and an exact sequence $0 \to N' \to N \to N'' \to 0$ such that $\mathrm{class}(F(N)) \in \Ext^1_{\cC_2}(F(N''), F(N'))$ coincides with $\mathrm{class}(V) \in \Ext^1_{\cC_2}(F(N'')), F(N'))$.

\medskip

We prove Claim 1. 
By the condition (vi) in \ref{C1C2} on semisimplicity, there is a morphism $W_wM =\gr^W_wM\to N'$  such that the composition $N'\to W_wM \to N'$ is the identity morphism. Let the exact sequence $0 \to N' \to N\to N'' \to 0$ be the pushforward of $0 \to W_wM \to U \to N'' \to 0$ under $W_wM \to N'$. Then this satisfies the condition in Claim 1.

\medskip

{\bf Claim 2.} $\mathrm{class}(U)\in \Ext^1_{\cC_1}(N'', W_wM)$ coincides with the image of $\mathrm{class}(N)\in \Ext^1_{\cC_1}(N'', N')$ under the homomorphism induced by $N'\to W_wM$.

\medskip
This follows from the injectivity of $\Ext^1_{\cC_1}(N'', W_wM) \to \Ext^1_{\cC_2}(F(N''), F(W_wM))$ (Lemma \ref{Claim0}) and the fact that $\mathrm{class}(F(U))\in \Ext^1_{\cC_2}(F(N''), F(W_wM))$ coincides with the image of $\mathrm{class}(V)\in \Ext^1_{\cC_2}(F(N''), F(N'))$ under the homomorphism  induced by the morphism $N' \to W_wM$. 

By Claim 2, there is a morphism $N\to M$ such that the diagram
$$\begin{matrix} 0 & \to & N' & \to &N & \to & N'' & \to & 0\\
&&\downarrow &&\downarrow &&\downarrow &&\\
0&\to& W_wM&\to& M &\to& M/W_wM&\to & 0
\end{matrix}$$
is commutative. This proves Proposition \ref{justi0}. 
\end{pf}

\begin{sbprop}\label{justi} Let the notation and the assumptions be as in $\ref{C1C2}$ (actually the condition {\rm (vi)} is not used for this proposition). Let $M$ and $N$ be objects of $\cC_1$. Then we have a bijection from $\Hom_{\cC_1}(M, N)$ to the set of pairs $(h, (h_w)_{w\in \Z})$, where $h$ is a morphism $F(M)\to F(N)$ and $h_w$ is a morphism $\gr^W_wM\to \gr^W_wN$ satisfying the following conditions {\rm (i)} and {\rm (ii)}.

{\rm (i)} $h$ sends $F(W_wM)$ to $F(W_wN)$ for any $w\in \Z$.

{\rm (ii)} For any $w\in \Z$, the morphism $F(\gr^W_wM)\to F(\gr^W_wN)$ induced by $h$ coincides with $F(h_w)$.
\end{sbprop}

\begin{pf} We first prove

\medskip

{\bf Claim 1.} Let $M$ and $N$ be objects of $\cC_1$. Let $w\in \Z$. Then we have a bijection from $\Hom_{\cC_1}(M, N)$ to the set of pairs $(a, b)$, where $a$ is a morphism $W_wM \to W_wN$ and $b$ is a morphism $M/W_wM\to N/W_wN$ satisfying the following condition $(*)$.

\medskip

$(*)$ The image of $\mathrm{class}(M)\in \Ext^1_{\cC_1}(M/W_wM, W_wM)$ in $\Ext^1_{\cC_1}(M/W_wM, W_wN)$ under the map induced by $a$ coincides with the image of $\mathrm{class}(N)\in \Ext^1_{\cC_1}(N/W_wN, W_wN)$ in $\Ext^1_{\cC_1}(M/W_wM, W_wN)$ under the map induced by $b$. 

\medskip
We prove Claim 1. Let $a: W_wM\to W_wN$ and $b: M/W_wM\to N/W_wN$ be morphisms. Let the exact sequence $0\to W_wN \to X \to M/W_wM\to 0$ be the pushforward of $0 \to W_wM \to M \to M/W_wM\to 0$ under $a$, and let the exact sequence $0\to W_wN \to Y \to M/W_wM \to 0$ be  the pullback of  $0\to W_wN \to N \to N/W_wN \to 0$ under $b$. Then 
the condition $(*)$ is that the extension classes of $X$ and $Y$ coincide. On the other hand, a morphism $M\to N$ which induces $a$ and $b$ corresponds bijectively to a morphism $X\to Y$ which induces identity morphisms of $W_wN$ and $M/W_wM$. By the first part of (2) of Lemma \ref{Claim0}, we have $\Hom_{\cC_1}(M/W_wM, W_wN)=0$. Hence we have the bijection in Claim 1. 
\medskip

We can prove similarly

\medskip
{\bf Claim 2.} Let $M$ and $N$ be objects of $\cC_1$. Let $w\in \Z$. Then we have a bijection from $\{h\in \Hom(F(M), F(N))\;|\; h(F(W_wM))\subset F(W_wN)\} $ to the set of pairs $(a, b)$, where $a$ is a morphism $F(W_wM) \to F(W_wN)$ and $b$ is a morphism $F(M/W_wM)\to F(N/W_wN)$ satisfying the following condition $(**)$.

\medskip

$(**)$ The image of $\mathrm{class}(F(M))\in \Ext^1_{\cC_2}(F(M/W_wM), F(W_wM))$ in $\Ext^1_{\cC_2}F((W_wM), F(W_wN))$ under the map induced by $a$ coincides with the image of $\mathrm{class}(F(N))\in  \Ext^1_{\cC_2}(F(N/W_wN), F(W_wN))$ in $\Ext^1_{\cC_2}(F(M/W_wM), F(W_wN))$ under the map induced by $b$. 

\medskip

Now we prove Proposition \ref{justi}. By downward induction on $w$, we may assume that there is a $w\in\Z$ such that $W_{w-1}M=M$, $W_{w-1}N=N$ and such that Proposition \ref{justi} is true if we replace $M$ and $N$ by  $M/W_wM$ and $N/W_wN$, respectively. By the Claim 1 and Claim 2, we have a commutative diagram with exact rows 
$$\begin{matrix} 0 &\to & A & \to & B & \to & C\\
&& \downarrow && \downarrow && \downarrow\\
0 &\to & A' & \to & B' & \to & C',
\end{matrix}$$
where 
$$A=\Hom_{\cC_1}(M, N), \quad A' =\Hom_{\cC_2,W}(F(M), F( N)) ,$$
$$B= \Hom_{\cC_1}(W_wM, W_wN) \times \Hom_{\cC_1}(M/W_wM, N/W_wN),$$
$$B'= \Hom_{\cC_2}(F(W_wM), F(W_wN)) \times \Hom_{\cC_2,W}(F(M/W_wM),  F(N/W_wN)),$$
$$C=\Ext^1_{\cC_1}(M/W_wM, W_wN), \quad C'=\Ext^1_{\cC_2}(F(M/W_wM), F(W_wN)).$$
Here $\Hom_{\cC_2,W}$ means the set of homomorphisms of $\cC_2$ which respect the filtrations $W$. The vertical arrows are injective by Lemma \ref{Claim0}. This proves 
$$A\overset{\cong}\to \{x\in A'\;|\; \text{the image of $x$ in $B'$ comes from $B$}\},$$ 
which proves Proposition \ref{justi} by downward induction on $w$. 
\end{pf}

\subsection{Main theorem}
\begin{sbpara}
Recall that the following (i) and (ii) are equivalent (Theorem \ref{thmss} (2)).

(i) In the category of log motives, homological equivalence (i.e.\ the trivial equivalence) coincides with the numerical equivalence.

(ii) The category of log motives is a semisimple abelian category.

\end{sbpara}
\begin{sbthm}\label{mthm}  {\rm (i)} and {\rm (ii)} are equivalent to the following {\rm (iii)}. 

{\rm (iii)} The category of log mixed motives is a  Tannakian category ({\rm \cite{Saavedra}, \cite{D3}}). 
\end{sbthm}

\begin{sbpara} We prove (ii) $\Rightarrow$ (iii). It is sufficient to prove that a morphism $f$ is an isomorphism if it induces an isomorphism $\cF\to \cF'$. By (ii), there is a morphism $h(X_w')\to h(X_w)$ which induces the inverse map $V'_{w,2}/V_{w,1}'\to V_{w,2}/V_{w,1}$. Thus the inverse map $\cF'\to \cF$ is a morphism of log mixed motives. 

We prove (iii) $\Rightarrow$ (i). Let $X$ be a projective vertical log smooth fs log scheme over $S$. Consider a morphism $f:h(X) \to h(X)$ which is numerically equivalent to $0$. We prove $f=0$. Let $V_1$ be the kernel of $f: h(X)_{\ell}\to h(X)_{\ell}$ and let $V_2=h(X)_{\ell}$. On the other hand, let $V'_1=0$ and $V'_2$ be the image of $f:h(X)_{\ell}\to h(X)_{\ell}$. Then $f$ induces an isomorphism $f: V_2/V_1\overset{\cong}\to V'_2/V'_1$. By (iii), there is a morphism $g: h(X) \to h(X)$ which induces the inverse map $V'_2/V'_1\to V_2/V_1$. Then $fg: h(X)_{\ell}\to h(X)_{\ell}$ is a projection to $V'_2$. Hence $\mathrm{Tr}(fg)=\dim(V_2')$. 
  Hence $\mathrm{Tr}(fg)=0$ implies $V_2'=0$ and hence $f=0$. 
\end{sbpara}

\section{Formulation with various realizations}
\label{realization} 
In Sections \ref{s:logmot} and \ref{s:logmm},  we considered $\ell$-adic realizations of log mixed motives fixing a prime number $\ell$. 

\subsection{Log motives and log mixed motives with many realizations}
\begin{sbpara}
  Let $\cR$ be the union of the set of all prime numbers and the set of three letters $\{B, D, H\}$: $B$ means Betti realization; $D$ means de Rham realization; $H$ means Hodge realization. 

  Let $S$ be an fs log scheme. 
  Let $R$ be a nonempty subset of $\cR$. 
  If a prime number $\ell$ is contained in $R$, assume that $S$ is over $\Z[1/\ell]$. 
  If $B\in R$, assume that $S$ is locally of finite type over $\C$. 
  If $D\in R$, assume that $S$ is log smooth over a field of characteristic $0$ or $S$ is a standard log point associated to a field of characteristic $0$. 
  If $H\in R$, assume that $S$ is log smooth over $\bC$ or $S$ is the standard log point associated to $\bC$. 
\end{sbpara}

\begin{sbpara}
  We define the categories 
$$\mathrm{LM}_R(S),\quad \mathrm{LMM}_R(S)$$
of log motives over $S$ and of log mixed motives over $S$, respectively, with respect to realizations in $R$. 

  The definition of $\mathrm{LM}_R(S)$ is similar to Section \ref{s:logmot}. 
  For a projective vertical log smooth fs log scheme $f\colon X \to S$ over $S$ and for $r\in \Z$, consider the symbol $h_R(X)(r)$. 

  When a prime $\ell$ belongs to $R$, let $$h_R(X)(r)_{\ell}:=\bigoplus_m\; H^m(X)_{\ell}(r).$$

  When $B \in R$, let $$h_R(X)(r)_{B}:=\bigoplus_m\; H^m(X)_{B}(r),\quad \mathrm{where}\quad H^m(X)_{B}= R^mf^{\log}_*\Q.$$ 
This is a locally constant sheaf of finite dimensional $\Q$-vector spaces on $S^{\log}$ (see Proposition \ref{prsmB}). 

  When $D \in R$, let $$h_R(X)(r)_{D}:=\bigoplus_m\; H^m(X)_{D}(r),\quad \mathrm{where}\quad H^m(X)_{D}= R^mf_{\mathrm{k\acute{e}t*}}
\omega_{X/S}^{\cdot,\mathrm{k\acute{e}t}}.$$
This is a locally free sheaf of $\cO_{\mathrm{k\acute{e}t}}$-modules of finite rank with a quasi-nilpotent integrable connection with log poles 
on $S_{\mathrm{k\acute{e}t}}$ (see Propositions \ref{prsmdR} (1) and \ref{prsmdR2} (1)). 

  When $H \in R$, let $$h_R(X)(r)_{H}:=\bigoplus_m\; H^m(X)_{H}(r),\quad \mathrm{where}\quad H^m(X)_{H}= R^mf_{\mathrm{k\acute{e}t*}}
\omega_{X/S}^{\cdot,\mathrm{k\acute{e}t}}$$
endowed with the natural log Hodge structures. 
  This is a log mixed Hodge structures on $S_{\mathrm{k\acute{e}t}}$ (see Propositions \ref{prsmdR} (2) and \ref{prsmdR2} (2)). 

  A morphism $h_R(X)(r) \to h_R(Y)(s)$ is defined as a family of morphisms between realizations for each element of $R$ satisfying,  
for any geometric standard log point $p$ over $S$, 
the pull-backed morphism is induced by a common element of gr of the $K$-group. 
  Note that we do not impose any comparison isomorphism between different realizations. 
  The rest is the same as in \ref{LM}, and we have the category $\mathrm{LM}_R(S)$.
  Here we use the Poincar\'e duality (Proposition \ref{PoinB}) and the K\"unneth formula (Proposition \ref{KunnB}) in log Betti cohomology, 
which implies the necessary corresponding theorems in log de Rham and log Hodge theory via log Riemann--Hilbert correspondence 
(\cite{IllusieKatoNakayama} {\sc Theorem} (6.2)). 
  We also use the Riemann--Roch theorems. 
\end{sbpara}

\begin{sbpara}
  The definition of $\mathrm{LMM}_R(S)$ is also similar to the case where $R$ consists of one prime. 
  We first define $\cC_{S,R}$ as follows. 

  First, for an $R$ consisting of one prime $\ell$, $\cC_{S,R}$ is $\cC_S$ in \ref{categoryC_S}.

  Second, for $R=\{B\}$, we define $\cC_{S,R}$ as the following category.

  Objects: $(\cF, W, (X_w)_{w\in \Z},  (V_{w,1})_{w\in \Z}, (V_{w,2})_{w\in \Z}, (\iota_w)_{w\in \Z})$. 

  Here $\cF$ is a locally constant sheaf of finite dimensional $\Q$-vector spaces on $S^{\log}_{\an}$. 
  The $W$ is an increasing filtration on $\cF$ by locally constant $\Q$-subsheaves. 
  The $X_w$ is a projective vertical log smooth fs log scheme over $S$. 
  For each $w\in \Z$, $V_{w,1}$ and $V_{w,2}$ are locally constant $\Q$-subsheaves of 
$\bigoplus_{r\in \Z} H^{w+2r}(X_w)_B(r)$ such that $V_{w,1}\subset V_{w,2}$. 
  The $\iota_w$ is an isomorphism $\gr^W_w\cF\cong V_{w,2}/V_{w,1}$. 

  A morphism $$(\cF, W, (X_w)_{w\in \Z}, (V_{w,1})_{w\in \Z}, (V_{w,2})_{w\in \Z}, (\iota_w)_{w\in \Z})\to (\cF', W', (X'_w)_{w\in \Z}, (V'_{w,1})_{w\in \Z}, (V'_{w,2})_{w\in \Z}, (\iota'_w)_{w\in \Z})$$ is a homomorphism of $\Q$-sheaves $\cF\to \cF'$ which respects the weight filtrations such that for each $w\in \Z$, the pullback of $\gr^W_w\cF\to \gr^{W'}_w\cF'$ to any geometric standard log point $s$ associated to $\bC$ over $S$ is induced from the sum of morphisms $h_{\{B\}}(X_w\times_S s)(r)\to h_{\{B\}}(X'_w \times_S s)(r')$ for various $r, r'$ 
which sends $V_{w,i}$ to $V'_{w,i}$ over $s$ for $i=1,2$. 

  Third, for $R=\{D\}$, we define $\cC_{S,R}$ as the following category.

  Objects: $(\cF, W, (X_w)_{w\in \Z},  (V_{w,1})_{w\in \Z}, (V_{w,2})_{w\in \Z}, (\iota_w)_{w\in \Z})$. 

  Here $\cF$ is a locally free $\cO_{S_{{\rm k\acute{e}t}}}$-modules of finite rank endowed with a quasi-nilpotent integrable connection with log poles.
  The $W$ is an increasing filtration on $\cF$ by locally free 
$\cO_{S_{{\rm k\acute{e}t}}}$-submodules with the compatible connections such that the graded quotients are also locally free.
  The $X_w$ is a projective vertical log smooth fs log scheme over $S$. 
  For each $w\in \Z$, $V_{w,1}$ and $V_{w,2}$ are locally free 
$\cO_{S_{{\rm k\acute{e}t}}}$-submodules with the compatible connections of 
$\bigoplus_{r\in \Z} H^{w+2r}(X_w)_D(r)$ such that $V_{w,1}\subset V_{w,2}$. 
  The $\iota_w$ is an isomorphism $\gr^W_w\cF\cong V_{w,2}/V_{w,1}$. 

  A morphism $$(\cF, W, (X_w)_{w\in \Z}, (V_{w,1})_{w\in \Z}, (V_{w,2})_{w\in \Z}, (\iota_w)_{w\in \Z})\to (\cF', W', (X'_w)_{w\in \Z}, (V'_{w,1})_{w\in \Z}, (V'_{w,2})_{w\in \Z}, (\iota'_w)_{w\in \Z})$$ is a homomorphism of 
$\cO_{S_{{\rm k\acute{e}t}}}$-modules $\cF\to \cF'$ which respects the weight filtrations such that for each $w\in \Z$, the pullback of $\gr^W_w\cF\to \gr^{W'}_w\cF'$ to any geometric standard log point $s$ over $S$ is induced from the sum of morphisms $h_{\{D\}}(X_w\times_S s)(r)\to h_{\{D\}}(X'_w \times_S s)(r')$ for various $r, r'$ 
which sends $V_{w,i}$ to $V'_{w,i}$ over $s$ for $i=1,2$. 

  Fourth, for $R=\{H\}$, we define $\cC_{S,R}$ as the following category.

  Objects: $(\cF, W, (X_w)_{w\in \Z},  (V_{w,1})_{w\in \Z}, (V_{w,2})_{w\in \Z}, (\iota_w)_{w\in \Z})$. 

  Here $(\cF, W)$ is a log mixed Hodge structure on $S_{\mathrm{k\acute{e}t}}$.   
  The $X_w$ is a projective vertical log smooth fs log scheme over $S$. 
  For each $w\in \Z$, $V_{w,1}$ and $V_{w,2}$ are sub-log Hodge structures of 
$\bigoplus_{r\in \Z} H^{w+2r}(X_w)_H(r)$ such that $V_{w,1}\subset V_{w,2}$. 
  The $\iota_w$ is an isomorphism $\gr^W_w\cF\cong V_{w,2}/V_{w,1}$. 

  A morphism $$(\cF, W, (X_w)_{w\in \Z}, (V_{w,1})_{w\in \Z}, (V_{w,2})_{w\in \Z}, (\iota_w)_{w\in \Z})\to (\cF', W', (X'_w)_{w\in \Z}, (V'_{w,1})_{w\in \Z}, (V'_{w,2})_{w\in \Z}, (\iota'_w)_{w\in \Z})$$ is a homomorphism of log mixed Hodge structures
$(\cF, W) \to (\cF', W')$ such that for each $w\in \Z$, the pullback of $\gr^W_w\cF\to \gr^{W'}_w\cF'$ to any standard log point $s$ associated to $\bC$ over $S$ is induced from the sum of morphisms $h_{\{H\}}(X_w\times_S s)(r)\to h_{\{H\}}(X'_w \times_S s)(r')$ for various $r, r'$ 
which sends $V_{w,i}$ to $V'_{w,i}$ over $s$ for $i=1,2$. 

  Lastly, for any $R$, we define $\cC_{S,R}$ as follows. 

  Objects: $(Y_{\rho})_{{\rho} \in R}$, where $Y_{\rho}$ is an object of $\cC_{S,\{ {\rho}\}}$, satisfying the condition that 
for any $w \in \bZ$, the $X_w$ of $Y_{\rho}$ is common. 
  
  A morphism $(Y_{\rho})_{{\rho} \in R} \to (Y'_{\rho})_{{\rho} \in R}$ is $(f_{\rho})_{{\rho}\in R}$, where 
$f_{\rho}\colon Y_{\rho} \to Y'_{\rho}$ is a morphism of $\cC_{S,\{ {\rho}\}}$, satisfying the condition that 
for any $w, r, r' \in \bZ$ and any $s \to S$, the element of gr of the $K$-group inducing the morphism $h_{{\rho}}(X_w \times_Ss)(r) \to
h_{{\rho}}(X'_w \times_Ss)(r')$ is common.

  Note that in this definition, we do not use any comparison isomorphism between different realizations. 
\end{sbpara}  

\begin{sbpara}  We define $\cC_{S,R}^{\mathrm{mot}}\subset \cC_{S,R}$ and $\mathrm{LMM}_R(S)$ imitating \ref{ss:LMM}. 
  Here the objects associated to standard objects for $B$, $D$, and $H$ are defined by virtue of Propositions \ref{prsmB}, 
\ref{prsmdR}, and \ref{prsmdR2}.
\end{sbpara}

\subsection{Conjectures and results}
  We state the conjecture that our categories $\mathrm{LM}_R(S)$ and $\mathrm{LMM}_R(S)$ are independent of the choices of the family $R$ of 
realizations. 
  We also state Tate conjecture and Hodge conjecture. 
  For the latter, we explain in Section \ref{example} that they hold in a simple case. 
  In there, we use the theories of log abelian varieties and log Jacobian varieties. 

\begin{sbconj} 
  Let $R'$ be a non-empty subset of $R$. Then the restriction of realizations give an equivalence of categories
$$\mathrm{LM}_R(S) \overset{\simeq}\to \mathrm{LM}_{R'}(S), \quad \mathrm{LMM}_R(S) \overset{\simeq}\to \mathrm{LMM}_{R'}(S).$$
\end{sbconj}

\begin{sbthm}
  The following {\rm (i)}--{\rm (iii)} are equivalent. 

{\rm (i)} In the category $\mathrm{LM}_R(S)$, homological equivalence (i.e.\ the trivial equivalence) coincides with the numerical equivalence.

{\rm (ii)} The category $\mathrm{LM}_R(S)$ is a semisimple abelian category. 

{\rm (iii)} The category $\mathrm{LMM}_R(S)$ is a Tannakian category. 
\end{sbthm}

\begin{pf}
  Similar to Theorems \ref{thmss} (2) and \ref{mthm}. 
\end{pf}

  For $\rho \in R$, we denote the realization for $\rho$ of $M\in \mathrm{LMM}(S)$ by $M_{\rho}$. 

\begin{sbconj}\label{Tconj} (Tate conjecture for log mixed motives.)
  Assume that $\ell$ is invertible over $S$. 
  Assume that either one of the following (i) and (ii) is satisfied. 

(i) $S$ is of finite type over some field which is finitely generated over the prime field.

(ii) $S$ is of finite type over $\Z$.

Then for any objects $M$ and $N$ of $\mathrm{LMM}_{\{\ell\}}(S)$, we have $$\Q_{\ell}\otimes \Hom(M, N) \overset{\cong}\to \Hom_W(M_{\ell}, N_{\ell}).$$ Here the right-hand-side denotes the set of homomorphisms of $\Q_{\ell}$-sheaves which respect the weight filtrations.
\end{sbconj}

\noindent 
{\it Remark.} If either the log structure of $S$ is trivial or $M$ and $N$ are pure, 
\lq\lq$W$'' on the right-hand-side in this conjecture can be eliminated (for the weight filtrations are automatically respected). 

\begin{sbconj}\label{2Tconj} (The second Tate conjecture.)
  Assume that $S$ is of finite type over $\Q$ and let  $M, N$ be objects of  $\mathrm{LMM}_{\{\ell, B\}}(S)$. 
  Then we have a bijection from 
$\Hom(M, N)$ to the set of all pairs $(a,b)$, where $a$ is a morphism $M_{\ell}\to N_{\ell}$, and 
$b$ is a homomorphism $M_B\to N_B$ defined on $(S\otimes \C)_{\an}^{\log}$, such that the pullback of $a$ on $(S\otimes \C)_{\an}^{\log}$ is induced from $b$ (cf.\ Proposition \ref{comparison}). 
\end{sbconj}

\begin{sbpara}
  The above second Tate conjecture follows from Tate conjecture. In fact, in
$\Q_{\ell} \otimes \Hom(M,N) \to \Q_{\ell} \otimes \{(a,b)\} \to \Hom(M_{\ell}, N_{\ell})$, the composition is an isomorphism if Tate conjecture is true and the second map is an injection. 
\end{sbpara}

\begin{sbconj}\label{Hconj}
(Hodge conjecture for log mixed motives.) 
  Assume that $S$ is log smooth over $\C$ or is the standard log point over $\bC$.
  Let $M$ and $N$ be objects of $\mathrm{LMM}_{\{H\}}(S)$. 
  Then we have
$$\Hom(M, N) \overset{\cong}\to \Hom(M_H, N_H).$$
\end{sbconj}

  By Proposition \ref{dRbasechange}, the conjecture \ref{Hconj} is reduced to the case where $S$ is the standard log point associated to $\C$.

\section{Examples}
\label{example}
\subsection{Log abelian varieties}
\label{la}
This \ref{la} and \ref{lj} are  preparations for \ref{curve} and \ref{opencurve}.
In this \ref{la}, we review the theory of log abelian varieties \cite{KKN2} and supply some results. 
  See \cite{Nakayamalasurvey} for a survey of the theory. 
We only consider log abelian varieties over a standard log point, for we need only this case in \ref{curve} and \ref{opencurve}. 

\begin{sbpara}\label{Lav1}
For an fs log scheme $S$, let $(\fs/S)$ be the category of fs log schemes over $S$, and let $(\fs/S)_{\mathrm{\acute{e}t}}$ be the site $(\fs/S)$ endowed with  the classical \'etale topology. A log abelian variety over $S$ is a sheaf of abelian groups on $(\fs/S)_{\mathrm{\acute{e}t}}$ satisfying certain conditions.
  If $s$ is the standard log point associated to a field $k$,  a log abelian variety over $s$ is described as in \ref{Lav2}--\ref{Lav6} below.   
\end{sbpara} 

\begin{sbpara}\label{Lav2}
Let ${\bf G}_{m,\log}$ be the sheaf $U \mapsto \Gamma(U, M^{\gp}_U)$ on $(\fs/s)_{\mathrm{\acute{e}t}}$. 

For a semiabelian variety $G$ over $k$ with the exact sequence $0\to T\to G \to B\to 0$, where $T$ is a torus over $k$ and $B$ is an abelian variety over $k$, let $G_{\log}$ be the pushout of $G \leftarrow T \to \cH om(X(T), {\bf G}_{m,\log})$ in the category of sheaves of abelian groups on $(\fs/s)_{\mathrm{\acute{e}t}}$. Here $X(T):= \cH om(T, {\bf G}_m)$ is the character group of $T$.  
We have $G\subset G_{\log}$.

Let $\cM_1$ be the category of systems $(\Gamma, G, h)$, where $\Gamma$ is a locally constant sheaf of free $\Z$-modules of finite rank on $(\fs/s)_{\mathrm{\acute{e}t}}$, $G$ is a semiabelian variety over $k$, and $h$ is a homomorphism $\Gamma\to G_{\log}$. 

An object of $\cM_1$ was called a log $1$-motif in \cite{KKN2}. 
\end{sbpara}

\begin{sbpara}\label{Lav3}  For an object $(\Gamma, G, h)$ of $\cM_1$  with $T$ the torus part of $G$, we have the $\Z$-bilinear paring 
$$\langle\;,\;\rangle: X(T) \times \Gamma \to \Z$$
(called the {\it monodromy pairing}) defined as follows. The map $h$ induces $\Gamma \to G_{\log} \to G_{\log}/G \cong T_{\log}/T$ and hence $X(T) \times \Gamma \to X(T) \times T_{\log}/T \to {\bf G}_{m, \log}/{\bf G}_m$. 
Since ${\bf G}_{m,\log}/{\bf G}_m$ restricted to the small \'etale site of $\Spec(k)$ is $\Z$, we have the above monodromy pairing. 

\end{sbpara}

\begin{sbpara}\label{Lav4}
  Let  $E=(\Gamma, G, h)$ be an object of $\cM_1$. The dual $E^*=(\Gamma^*, G^*, h^*)$ of $E$ is an object of $\cM_1$ defined as in \cite{KKN2}.  We have $\Gamma^*= X(T)$, the torus part $T^*$ of $G^*$ is $\cH om(\Gamma, {\bf G}_m)$, and the abelian variety $G^*/T^*$  is the dual abelian variety $B^*$ of $B=G/T$. 

Let  $E=(\Gamma, G, h)$ be an object of $\cM_1$. 

A {\it polarization} on $E$ is a homomorphism $p: E\to E^*$ satisfying the following conditions (i)--(iv). 

\medskip
(i) The homomorphism $B\to B^*$ induced by $p$ is a polarization of the abelian variety $B$. 

(ii) The homomorphism $\Gamma \otimes \Q\to \Gamma^*\otimes \Q$ induced by $p$ is an isomorphism. 

(iii) The pairing $\Gamma \times \Gamma \to \Z\;;\; (a, b)\mapsto \langle p(a), b \rangle$ is a positive definite symmetric bilinear form, where $\langle\;,\;\rangle$ denotes the monodromy pairing (\ref{Lav3}) and $p$ denotes the homomorphism $\Gamma \to \Gamma^*=X(T)$ induced by $p$.  

(iv) The  homomorphism $T_{\log}\to (T^*)_{\log}$  induced by $p$ comes from $T\to T^*=\cH om(\Gamma, {\bf G}_m)$ which is dual to the homomorphism $\Gamma \to \Gamma^*=X(T)$ induced by $p$.

\medskip

Let $\cM_0$ be the full subcategory of $\cM_1$ consisting of objects which have polarizations after the base change to $\overline k$. 
\end{sbpara}

\begin{sbpara}\label{Lav6}  For an object $(\Gamma, G, h)$ of $\cM_1$, we have a subgroup sheaf  $G_{\log}^{(\Gamma)}$
of $G_{\log}$ containing $G$ and $h(\Gamma)$   defined as in \cite{KKN2}.

A {\it log abelian variety over $s$} is a sheaf of abelian groups $A$ on $(\fs/s)_{\mathrm{\acute{e}t}}$ such that $A= G_{\log}^{(\Gamma)}/h(\Gamma)$ for some object  $(\Gamma, G, h)$ of $\cM_0$. Let $\mathrm{LAV}(s)$ be the category of log abelian varieties over $s$. 
We have an equivalence of categories  $$\cM_0 \overset{\sim}\to \mathrm{LAV}(s)\;;\; (\Gamma, G, h) \mapsto G_{\log}^{(\Gamma)}/h(\Gamma)$$ 
by \cite{KKN2} {\sc Theorem} 3.4 (cf.\ \cite{KKN2} {\sc Proposition} 4.5 and \cite{KKN2} {\sc Theorem} 4.6 (2)). 
\end{sbpara}
\begin{sbpara} Let $E$ be an object of $\cM_0$ and let $A$ be the corresponding log abelian variety. Then the  log abelian variety $A^*$ corresponding to the dual $E^*$ of $E$ is called the {\it dual log abelian variety} of $A$. We have an embedding $A^*\subset \cE xt^1(A, {\bf G}_{m,\log})$. A polarization of $A$ gives a homomorphism $A\to A^*$. 
\end{sbpara}

\begin{sbpara}
For an additive category $\cC$, let $\cC \otimes \Q$ be the following category. 
Objects of $\cC \otimes \Q$  are the same as those of $\cC$. For objects $E, E'$ of $\cC$, $\Hom_{\cC \otimes \Q}(E, E')=\Hom_{\cC}(E, E')\otimes \Q$. 
\end{sbpara} 

\begin{sbpara} 
  The category $\cM_1 \otimes \Q$ is an abelian category as is seen easily.  
  $\cM_0 \otimes\Q$ is stable in $\cM_1\otimes \Q$ under taking kernels, cokernels, and direct sums (cf.\ \cite{Zhao}), 
and hence, it is an abelian category. Hence  $\mathrm{LAV}(s)\otimes \Q$  is an abelian category. 
\end{sbpara}

\begin{sbpara} Let $A$ be a log abelian variety over $s$ corresponding to an object $(\Gamma, G, h)$ of $\cM_0$. Let $\ell$ be a prime number which is different from the characteristic of $k$. Then the $\ell$-adic Tate module $T_{\ell}A$ is defined in the natural way as a smooth $\Z_{\ell}$-sheaf on the log \'etale site of $s$ (cf.\ \cite{KKN4} 18.9). Let $V_{\ell}A=\Q_{\ell}\otimes T_{\ell}A$. 

We have an exact sequence $0\to T_{\ell}G \to T_{\ell}A\to \Gamma \otimes \Z_{\ell}\to 0$ (cf.\ \cite{KKN4} 18.10).

We have $T_{\ell}(A^*)= \cH om(T_{\ell}A, \Z_{\ell}(1))$. 

If $T$ is the torus part of $G$,  the monodromy operator $\cN :T_{\ell}A \to T_{\ell}A(-1)$ coincides with the composition $T_{\ell}A \to \Gamma \otimes \Z_{\ell} \to T_{\ell}T(-1)\to T_{\ell}A(-1)$, where the second arrow $\Gamma \otimes \Z_{\ell} \to T_{\ell}T(-1) = \Hom(X(T), \Z_{\ell})$ is the map induced by the monodromy pairing $\langle\;,\;\rangle: X(T) \times \Gamma \to \Z$ (\ref{Lav3}). 

\end{sbpara}

\begin{sbpara} 
Let $A$ be a polarizable log abelian variety over $s$. 
  Fix a polarization $p: A \to A^*$.
  Then $p$ is an isomorphism in $\mathrm{LAV}(s) \otimes \Q$. For $f \in \mathrm{End}_{\mathrm{LAV}(s) \otimes\Q}(A)$, let
 $f^{\sharp}:=p^{-1}f^*p\in \mathrm{End}_{\mathrm{LAV}(s)\otimes\Q}(A)$, where $f^*: A^*\to A^*$ is the dual of $f$. 
\end{sbpara}

\begin{sbprop} Let $A$ and $p$ be as above and let $f\in  \mathrm{End}_{\mathrm{LAV}(s)\otimes\Q}(A)$,  $f\neq 0$. Then 
$\Tr(ff^{\sharp})>0$. Here $\Tr$ is the trace of the induced $\Q_{\ell}$-linear map $V_{\ell}A \to V_{\ell}A$. 
 
 \end{sbprop}

\begin{pf} Let $E=(\Gamma, G, h)$ be an object of $\cM_0$ corresponding to $A$, let $T$ be the torus part of $G$, and let $B=G/T$ be the quotient abelian variety of $G$. Let $f_0, f_0^{\sharp} : \Gamma \otimes \Q_{\ell}\to \Gamma \otimes\Q_{\ell}$, $f_1, f_1^{\sharp}: V_{\ell}B \to V_{\ell}B$, and $f_2, f_2^{\sharp}: V_{\ell}T\to V_{\ell}T$ be the map induced by $f$, $f^{\sharp}$, respectively. Then 
$$\Tr(ff^{\sharp})= \sum_{i=0}^2 \; \Tr(f_if_i^{\sharp}).$$
By the usual theory of abelian varieties, $\Tr(f_1f_1^{\sharp})\geq 0$ and it is non-zero if $f_1\neq 0$. $\Tr(f_0 f_0^{\sharp})\geq 0$ and this is non-zero if $f_0 \neq 0$, for we have a positive definite symmetric form. We have $\Tr(f_2f_2^{\sharp})\geq 0$ and it is non-zero if $f_2\neq 0$ by duality. Hence  $\Tr(f f^{\sharp})\geq 0$ and this is non-zero unless $f_0=f_1=f_2=0$.  If $f_0=f_1=f_2=0$, $f=0$ because any homomorphism $B\to T$ is zero.
\end{pf}

\begin{sbcor}
\label{LAss}
The category $\mathrm{LAV}(s)\otimes \Q$ is semisimple. 
\end{sbcor}

\begin{pf}
  This is deduced from the above proposition by the arguments in \ref{semisimple}.
\end{pf}

\begin{sbpara} Let $A$ be a log abelian variety over $s$. Assume $k=\C$. Then 
 we have  the polarizable log Hodge structure 
  over $s$ of weight $-1$ corresponding to $A$  (\cite{KKN1}), which we denote by $H_1(A)_H$. The underlying 
 locally constant sheaf 
 of finite-dimensional $\Q$-vector spaces
  on the topological space $s_{\an}^{\log}$ (which is homeomorphic to a circle $S^1$) will be denoted by $H_1(A)_B$. If $(\Gamma, G, h)$ denotes the object of $\cM_0$ corresponding to $A$, we have an exact sequence
  $$0\to \cH_1(G, \Z) \to H_1(A)_B \to \Gamma \to 0.$$
  \end{sbpara}

\begin{sbprop} Let $A_1$ and $A_2$ be log abelian varieties over $s$.

$(1)$ If $k$ is finitely generated over a prime field, we have $\Z_{\ell}\otimes_\Z \Hom(A_1, A_2)\overset{\cong}\to \Hom(T_{\ell}A_1, T_{\ell}A_2)$. 

$(2)$ If $k$ is a subfield of $\C$ which is finitely generated  over $\Q$, we have a bijection from $\Hom(A_1, A_2)$ to the set of pairs $(a, b)$, where $a$ is a homomorphism $T_{\ell}A_1 \to T_{\ell}A_2$ and $b$ is a homomorphism $H_1(A_1)_B \to H_1(A_2)_B$ on $(s \otimes_k \C)^{\log}$ such that the pullback of $a$ on $(s \otimes_k \C)^{\log}$ is induced by $b$.

$(3)$ If $k=\C$, $\Hom(A_1, A_2) \overset{\cong}\to \Hom(H_1(A_1)_H, H_1(A_2)_H)$. 
\end{sbprop}

\begin{pf} 
 For an object $E=(\Gamma, G, h)$ of $\cM_1$, define the filtration $W$ on $E$ by $W_wE=E$ for $w\geq 0$, $W_{-1}E= (0, G, 0)$, $W_{-2}E=(0, T, 0)$ with $T$ the torus part of $G$, and $W_wE=0$ for $w\leq -3$. Then $\gr^W_0E=(\Gamma,0,0)$, $\gr^W_{-1}E=(0, B,0)$, where $B$ is the abelian variety $G/T$, $\gr^W_{-2}E=(0,T,0)$, and $\gr^W_wE=0$ for $w\neq 0,-1,-2$. Let $\cC_1= \cM_1\otimes\Q$ and let $\cC_2$ be the category of smooth $\Q_{\ell}$-sheaves on the  log \'etale site of $s$. Then (1) and (2) follow from the Tate conjecture on homomorphisms  of abelian varieties proved by Faltings (\cite{Faltings}) and from the injectivity of $G(k) \otimes\Q \to H^1(k, V_{\ell}G)$ for a semiabelian variety $G$ over $k$,
by the method of \ref{justification}.

(3) follows from \cite{KKN1}. 
\end{pf}

\subsection{Log Jacobian varieties}
\label{lj}
We review the theory of log Jacobian varieties of log curves over a standard log point in \cite{Kajiwara}, and supply some results. 
  In this subsection and the next, we omit some details of proofs, which will be treated in a forthcoming paper. 

\begin{sbpara} 
  Let $s$ be the standard log point associated to a field $k$. Let $X$ be a projective vertical log smooth connected curve over $s$ which is semistable, whose double points are rational and whose components are geometrically irreducible. 

Then we have a log abelian variety over $s$ associated to $X$ called the {\it log Jacobian variety} of $X$. 
  We will denote it by $J$. 

  This $J$ is essentially constructed by Kajiwara in \cite{Kajiwara}. 
  We explain his construction below in \ref{lj_construction}. 

  This $J$ has the following properties \ref{degree_map} and \ref{lj_1-motif}. 
\end{sbpara}

\begin{sbpara}
\label{degree_map} 
  Let $\cH^1(X, M^{\gp})$  be the sheafification of the presheaf $U \mapsto H^1(X \times_s U, M^{\gp})$  on $(\fs/s)_{\mathrm{\acute{e}t}}$. We have a degree map $\cH^1(X, M^{\gp})\to \Z$. Let $\cH^1(X, M^{\gp})^0\subset \cH^1(X, M^{\gp})$ be the kernel of the degree map. Then $J$ is a subgroup sheaf of 
 $\cH^1(X,  M^{\gp})^0$. 
\end{sbpara}
 
\begin{sbpara}
\label{lj_1-motif} 
  Let $E=(\Gamma, G, h)$ be the object of $\cM_0$ corresponding to $J$, let $T$ be the torus part of $G$, and let $B=G/T$ be the quotient abelian variety of $G$. Then $\Gamma$, $T$, $B$ are described as follows. 
 
  Let $\Gamma$ be the first homology group of the graph of $X$ as usual, that is, $\Gamma= \Ker(\bigoplus_{I_1} \Z \to \bigoplus_{I_0} \Z)$, where $I_0$ is the set of generic points of $X$, and  $I_1$ is the set of singular points of $X$. ${\cH} om\,(T,{\Bbb G}_m)= \Hom(\Gamma, \Z)$. $B= \prod_{\nu\in I_0} J_{D(\nu)}$, where $D(\nu)$ is the closure of $\nu$ in $X$ which is a projective smooth curve over $k$ and $J_{D(\nu)}$ is the Jacobian variety of $D(\nu)$. We have a canonical  isomorphism $J\cong J^*$ (cf.\ \ref{Al4}) which induces the evident isomorphisms $\Gamma \cong \Gamma^*$, $T\cong T^*$ and $B\cong B^*$. 
\end{sbpara}

\begin{sbpara}
\label{lj_construction}
  We explain the construction of $J$, which is essentially due to Kajiwara. 
  For simplicity, we assume that $k$ is algebraically closed. 
  By \cite{Kajiwara} (2.18), we have a commutative diagram 
$$\begin{CD}
@. 0  @. 0   @. 0 @.@. \\
@. @VVV @VVV  @VVV @. \\
@.\Gamma @. G_{\log} @. \cH^1(X,M^{\gp}) @. \\
@. @VVV @VVV  @| @. \\
0@>>> \bigoplus_{I_1} \bZ @>>> P^{\log}_{X/s} @>>> \cH^1(X,M^{\gp}) @>>> 0 \\
@. @VVV @VVV @VVV \\
0@>>> \bigoplus_{I_0} \bZ @= \bigoplus_{I_0} \bZ @>>>  0 \\
@. @VVV @. @. \\
@. \bZ @.@.@.
\end{CD}$$
with exact rows and columns, where $G=\Ker(\cH^1(X,\Gm)\to \bigoplus_{I_0}\bZ)$, and $P^{\log}_{X/s}$ is defined in \cite{Kajiwara}. 
  This diagram yields a log $1$-motif $(\Gamma, G, h\colon \Gamma \to G_{\log})$ and the degree map $\cH^1(X,M^{\gp}) \to \bZ$ whose kernel  $\cH^1(X, M^{\gp})^0 \cong G_{\log}/h(\Gamma)$ contains $G_{\log}^{(\Gamma)}/h(\Gamma)$. 
  The last sheaf is $J$. 
\end{sbpara}

\begin{sbpara}\label{Al2}
  Let $Y:= X \times_s X$. We have a $M^{\gp}_Y$-torsor on $Y$ called the {\it Poincar\'e torsor}, defined as follows.  

  Let $I$ be the ideal of $\cO_Y$ which defines the diagonal $X$ in $Y$. 
  Let $R$ be the sheaf of rings on $Y$ defined as the  $\cO_Y$-algebra generated over $\cO_Y$ by $f^{-1}$ for non-zero divisors $f$ in $\cO_Y$. 
  We have $M_Y^{\gp}\subset R^\times$. 
  There is a unique global section $t$ of $R^\times/M_Y^{\gp}$ having the following property: Locally on $Y$, $I$ is generated by local sections of  $R^\times \cap I$ which belong to $t$. 
  In fact, let $\pi$ be a generator of the standard log point. 
  At a singular point $x$ of $X$, 
let $f_1, g_1$ be generators of the log of the left $X$ in $X \times_sX$ around $x$ such that $f_1g_1 = \pi$, 
and let $f_2, g_2$ be the copies of them for the right $X$ in $X \times_sX$. 
  Then, $I$ is generated by $f_1-f_2$ and $g_1-g_2$. 
  We have $g_1-g_2= (-\pi f_1^{-1}f_2^{-1})(f_1-f_2)$ in $R^\times$ and $-\pi f_1^{-1}f_2^{-1}\in M^{\gp}_Y$. 
  The desired $t$ is defined as the class of $f_1-f_2$ which is also the class of $g_1-g_2$.

  Let the Poincar\'e torsor be 
the inverse image of $t^{-1}$ in $R^\times$ under $R^\times \to R^\times/M^{\gp}_Y$. This is an $M^{\gp}_Y$-torsor. 

  If $X$ comes from a projective smooth curve over $k$, this Poincar\'e torsor comes from the usual Poincar\'e ${\bf G}_m$-torsor. 
\end{sbpara}

\begin{sbpara}\label{Al3}
  We have a morphism $\varphi: X\to \cH^1(X, M^{\gp})$ which sends $x$ to the pullback of the Poincar\'e torsor (\ref{Al2}) 
with respect to $X \to X\times X; y \mapsto (x,y)$. 

  If $b$ is a morphism $s\to X$ over $s$, we have a canonical morphism 
$$\varphi_b:X\to J\subset \cH^1(X, M^{\gp})\;;\;x \mapsto \varphi(x)-\varphi(b)$$
called the {\it log Albanese mapping} associated to $b$. 
\end{sbpara}

\begin{sbpara}\label{Al4} (Self-duality of the log Jacobian.)
  Let $b$ and $\varphi_b$ be as above. Then the pulling back via $\varphi_b$ gives an
isomorphism $$\cE xt^1(J, {\bf G}_{m,\log})\overset{\cong}\to \cH^1(X, M^{\gp})^0,$$ which is independent of the choice of $b$. 
Hence the subgroup sheaf $J$ of $ \cH^1(X, M^{\gp})^0$ is regarded as a subgroup sheaf of $\cE xt^1(J, {\bf G}_{m,\log})$.   
  Via this inclusion $J \subset \cE xt^1(J, {\bf G}_{m,\log})$, $J$ is identified with the dual log abelian variety $J^*$ of $J$. 
  Since this isomorphism $J\cong J^*$ does not depend on $b$,  it is defined canonically even if there is no $b$. 
\end{sbpara}

\begin{sbprop}\label{Al5}
  Let $b:s \to X$ be a morphism over $s$.  Let $A$ be any log abelian variety over $s$. Then the map $$\Hom(J, A)\to \mathrm{Mor}(X, A)\;;\; h \mapsto h \circ \varphi_b$$ is bijective. 
\end{sbprop}
  
\begin{pf} 
  The inverse map is given as follows. Let 
  $f: X \to A$ be a morphism. Then we have  $A^*\to \cE xt^1(A, {\bf G}_{m,\log}) \to \cH^1(X, M^{\gp})^0$, where the second arrow is the pullback by $f$. This induces $A^*\to J$. Taking the dual log abelian varieties, we have $J\to A$. 
\end{pf}
 
\begin{sbpara}
\label{TateJ}
 Let $\ell$ be a prime number which is invertible in $k$. Then we have  canonical isomorphisms
 $$V_{\ell}J \cong H^1(X)_{\ell}(1)\cong  \cH om(H^1(X)_{\ell}, \Q_{\ell}).$$
\end{sbpara}

\subsection{Examples I}
\label{curve}
This subsection Examples I is for the pure case. The next subsection Examples II is for the mixed case. 

The following is a part of Conjecture \ref{conjoplus}.

\begin{sbprop}
\label{curveoplus}
  Let $X$ be a projective vertical log smooth curve over an fs log scheme $S$. Then $h(X)= h^0(X) \oplus h^1(X) \oplus h^2(X)$.
\end{sbprop}

\begin{pf}
  Using the $\ell$-adic weights, we reduce to the case where $S$ is a geometric standard log point. 
  Then we see that the desired morphisms come from the $K$-group by the existence of a section. 
\end{pf}

\begin{sbprop}\label{H2} 
  Assume that $S$ is %
the standard log point over $\bC$. 
  Let $X$ be a connected projective semistable %
curve over $S$.
Then the Hodge conjecture $\ref{Hconj}$ for $\Hom(\Q, h^2(X)(1))$ is true. 
\end{sbprop}

\begin{pf} %
Assume that we are given a homomorphism $h:\Q\to H^2(X)_B(1)$. By invariant cycle theorem, this comes from the classical Betti cohomology $H^2(X_{\an}, \Q(1))$. Since  $h(1)$ belongs to $\mathrm{Fil}^1H^2(X)_H$, it vanishes in $H^2(X, \cO_X)$. Hence it comes from the kernel of $H^2(X_{\an}, \Q(1))\to H^2(X, \cO_X)$. By the exponential sequence $0\to \Z(1) \to \cO_{X_{\an}} \to \cO_{X_{\an}}^\times\to 0$, it comes from
$\mathrm{Pic}(X) \otimes \Q$. 
\end{pf}

  The next will be proved in a forthcoming paper. 

\begin{sbprop}\label{CandJ}  Let $s$ be  a geometric standard log point of characteristic $\neq \ell$.  For $i=1,2$, let $X_i$ be a projective vertical log smooth curve over $s$ which is semistable, and let $J_i$ be the log Jacobian variety of $X_i$. For a homomorphism $h: H^1(X_1)_{\ell}\to H^1(X_2)_{\ell}$, the following two conditions {\rm{(i)}} and {\rm{(ii)}} are equivalent. 

{\rm{(i)}} $h$ is a morphism $H^1(X_1)\to H^1(X_2)$ of log motives over $s$. 

{\rm{(ii)}} $h$ comes from a morphism $J_1\to J_2$ in $\mathrm{LAV}(s)$ (via the isomorphisms $H^1(X_i)_{\ell}(1) \cong V_{\ell}J_i$ in $\ref{TateJ}$). 
\end{sbprop}

\begin{sbprop}\label{propE1}
 Let $X$ and $Y$ be projective vertical log smooth curves over an fs log scheme $S$ whose geometric fibers are connected. 
 
$(1)$  Assume that $S$ is %
the standard log point over $\bC$ and that $X$ and $Y$ are semistable over $S$. 
  Then
the Hodge conjecture $\ref{Hconj}$ for $\Hom(h(X), h(Y))$ is true. 

$(2)$ Assume that $S$ is of finite type over $\Q$. Then the second Tate conjecture 
$\ref{2Tconj}$ for 
$\Hom(h(X), h(Y))$ is true. 

$(3)$ Assume that $S$ is a standard log point associated to a finitely generated field over a prime field whose characteristic is different from a prime number $\ell$. 
  Then the Tate conjecture $\ref{Tconj}$ for 
$\Hom(h(X), h(Y))$ is true. 

$(4)$ For $f, g \in \Hom(h(X), h(Y))$, if $f$ and $g$ are numerically equivalent, then $f=g$.

$(5)$ The endomorphism ring of $h(X)$ is a finite-dimensional semisimple algebra over $\Q$. 
\end{sbprop}

\begin{pf} By \ref{TS}, %
Proposition \ref{prsmbcB} and 
$\ell$-adic log proper base change theorem \cite{KajiwaraNakayama} 
{\sc Proposition} 5.1 (cf.\ \cite{KajiwaraNakayama} {\sc Remark} 5.1.1), we may assume that $S$ is a standard log point and $X$ and $Y$ are semistable and that their double points are rational and their components are geometrically irreducible. 
  Let $J$ and $J'$ be the log Jacobian variety of $X$ and $Y$, respectively. 
  By Propositions \ref{curveoplus}, \ref{CandJ}, and the method of \ref{justification}, we can identify $\Hom(h(X), h(Y))$ with 
$\Hom_{\mathrm{LAV}(S)\otimes \bQ}(J, J')$. 
  Then we reduce to the results in \ref{la}.
\end{pf}

\subsection{Examples II}
\label{opencurve}
\begin{sbpara}\label{5.5.1}

Let $X$ be a projective vertical  log smooth curve over an fs log scheme $S$. Let $n\geq 1$ and $s_1, \dots, s_n: S\to X$ be strict morphisms over $S$ such that 
$s_i(S) \cap s_j(S)=\emptyset$ if $i\neq j$.
   Let $D:=\bigcup_{i=1}^n s_i(S)$ and let $U:= X\smallsetminus D$.

We will denote the log mixed motive corresponding to the standard object associated to $(U, X, D, 1)$ over $S$ by $H^1(U)$. 

Let $\Gamma = \mathrm{Ker}(\mathrm{sum}: \Z^n \to \Z)$. We have $W_0M=M$, $W_{-2}M=0$, $W_{-1}M= H^1(X)$, $\gr^W_0 M= \Gamma \otimes \Q(-1)$, where $M=H^1(U)$. 

The spectral sequence as 
 in \ref{mix8} for each realization degenerates at $E_2$. 
\end{sbpara}

\begin{sbpara}\label{5.5.2}
   Let the notation be as in \ref{5.5.1}.

  If $S$ is over $\bZ[1/\ell]$, we have an exact sequence

(1) $0 \to H^1(X)_{\ell} \to H^1(U)_{\ell} \to \Gamma \otimes \Q_{\ell}(-1) \to 0$

\noindent of $\Q_{\ell}$-sheaves. 

  If $S$ is either log smooth over $\bC$ or the standard log point associated to $\bC$, we have an exact sequence 

(2) $0 \to H^1(X)_H \to H^1(U)_H \to \Gamma \otimes \Q(-1) \to 0$

\noindent of log mixed Hodge structures over $S$. 

Assume that $S$ is a standard log point associated to a field $k$, and assume that $X$ is connected and semistable and that their double points are rational and their components are geometrically irreducible. 
  Let $J$ be the log Jacobian variety of $X$. Then $(s_i)_i$ induces a homomorphism $\psi:=\varphi\circ (s_i)_i: \Gamma \to J$ by the log Albanese mapping $\varphi$ (\ref{Al3}). 

Note that for any log abelian variety $A$ over $S$, we have a canonical homomorphism

(3) $A(S) \otimes \Q \to \mathrm{Ext}^1(\Q_{\ell}, V_{\ell}A)$

\noindent by Kummer theory, which is injective if $k$ is finitely generated over a prime field. If $k=\C$, we have also a canonical injective map

(4) $A(S)  \to \mathrm{Ext}^1(\Z, H_1(A)_H)$. 

We have:

(5) Under the homomorphism $\Hom(\Gamma, J) \to \mathrm{Ext}^1(\Gamma \otimes \Q_{\ell}, V_{\ell}J)$ induced by (3) (applied to the log abelian variety $A=\cH om(\Gamma, J)$), the extension class of (1) coincides with the image of $\psi: \Gamma \to J$. 

(6) If $k=\bC$, under the homomorphism $\Hom(\Gamma, J) \to \mathrm{Ext}^1(\Gamma, H_1(J)_H)$ induced by (4), the extension class of (2) coincides with the image of $\psi: \Gamma \to J$. 

\end{sbpara}

\begin{sbprop}\label{propE2}
 Let $U_1$, $U_2$ be objects as  $U$ in $\ref{5.5.1}$.  
 
$(1)$  Assume that $S$ is the standard log point over $\bC$ and that $X_1$ and $X_2$ are connected and semistable. 
  Then
the Hodge conjecture $\ref{Hconj}$ for $\Hom(H^1(U_1), H^1(U_2))$ is true. 

$(2)$ Assume that $S$ is of finite type over $\Q$. Then the second Tate conjecture 
$\ref{2Tconj}$ for 
$\Hom(H^1(U_1), H^1(U_2))$ is true. 

$(3)$ Assume that $S$ is the standard log point associated to a finitely generated field over a prime field whose characteristic is different from a prime number $\ell$. 
  Then the Tate conjecture $\ref{Tconj}$ for 
$\Hom(H^1(U_1), H^1(U_2))$ is true. 
\end{sbprop}

\begin{pf} Similarly as in Proposition \ref{propE1}, we may assume that $S$ is a standard log point and $X_i$ are connected and semistable and that their double points are rational and their components are geometrically irreducible.  For $i=1, 2$, let  $J_i$ be the log Jacobian variety of $X_i$. By (5) in \ref{5.5.2}, by the injectivity of the map (3) in \ref{5.5.2},  and by Proposition \ref{CandJ}, the method of \ref{justification} shows  that 

$(*)$ the set of morphisms 
$H^1(U_1)\to H^1(U_2)$ is identified with the set of pairs $(a, b)$, where $a$ is a homomorphism $\Gamma_1 \otimes \Q \to \Gamma_2\otimes \Q$ and $b$ is a  morphism $J_1 \to J_2$ in $\mathrm{LAV}(s)\otimes \Q$ such that $\psi_2\circ a = b \circ \psi_1$.

 Hence by (5) in \ref{5.5.2}, by the injectivity of the map (3) in \ref{5.5.2}, and by this $(*)$, the method of \ref{justification} proves (2) and (3). 
  Similarly, by (6) in \ref{5.5.2}, by the injectivity of the map (4) in \ref{5.5.2}, and by $(*)$, the method of \ref{justification} proves (1). 
\end{pf}

\bigskip

\noindent {\rm Tetsushi ITO
\\
Department of Mathematics
\\
Kyoto University
\\
Kitashirakawa, Kyoto 606-8502, Japan}
\\
{\tt tetsushi@math.kyoto-u.ac.jp} 

\bigskip

\noindent {\rm Kazuya KATO
\\
Department of Mathematics
\\
University of Chicago
\\
Chicago, Illinois, 60637, USA}
\\
{\tt kkato@math.uchicago.edu}

\bigskip

\noindent {\rm Chikara NAKAYAMA
\\
Department of Economics 
\\
Hitotsubashi University 
\\
2-1 Naka, Kunitachi, Tokyo 186-8601, Japan}
\\
{\tt c.nakayama@r.hit-u.ac.jp}

\bigskip

\noindent
{\rm Sampei USUI
\\
Graduate School of Science
\\
Osaka University
\\
Toyonaka, Osaka, 560-0043, Japan}
\\
{\tt usui@math.sci.osaka-u.ac.jp}
\end{document}